\documentclass[a4paper,11pt,oneside,fleqn]{article}

\usepackage[utf8]{inputenc}
\usepackage[ngerman,english]{babel}
\usepackage[T1]{fontenc}
\usepackage{lmodern}
\usepackage{standalone}
\usepackage{mathtools}
\usepackage{amsfonts}
\usepackage{amssymb}
\usepackage{hyperref}
\usepackage{enumerate}
\usepackage{amscd}
\usepackage{graphicx}
\usepackage{pgfplots}
\usepackage{caption}
\usepackage{subcaption}
\usepackage{extarrows}
\usepackage{faktor}
\usepackage[vcentermath]{youngtab}
\usepackage{marvosym}
\usepackage{mathtools}
\usepackage{color}
\usepackage{colortbl}
\usepackage{makeidx}
\usepackage[all]{xy}
\usepackage{amsthm}
\usepackage{paralist}
\usepackage{enumitem}
\usepackage{tikz,pgf}
\usepackage{tikz-cd}
\usetikzlibrary{patterns}
\usetikzlibrary{matrix, arrows}

\usepackage{soul}

\hypersetup{colorlinks, linkcolor=blue, citecolor=blue, urlcolor=blue}

\usepackage{stmaryrd}

\usepackage{multicol}

\usepackage[sort,numbers]{natbib}
\allowdisplaybreaks

\makeatletter
\long\def\@makecaption#1#2{%
  \vskip\abovecaptionskip\footnotesize
  \sbox\@tempboxa{#1. #2}%
  \ifdim \wd\@tempboxa >\hsize
    #1. #2\par
  \else
    \global \@minipagefalse
    \hb@xt@\hsize{\hfil\box\@tempboxa\hfil}%
  \fi
  \vskip\belowcaptionskip}
\makeatother

\newcommand{\todo}[1][\null]{\ensuremath{\clubsuit}}

\newcommand{\noprint}[1]{}
\newcommand{\checked}[1][\null]{\ensuremath{\boldsymbol{\surd}}}
\newcommand{\R}{\mathbb{R}}
\newcommand{\M}{\mathbb{M}}
\newcommand{\h}{H}
\newcommand{\m}{\mathbf{m}}

\newcommand{\mybrL}{\llbracket}
\newcommand{\mybrR}{\rrbracket}

\renewcommand{\u}{\mathbf{u}}

\renewcommand{\d}{\mathbf{d}}

\newcommand{\eps}{\varepsilon}

\newcommand{\C}{\mathcal{C}}
\DeclareMathOperator\Div{Div}
\DeclareMathOperator\Curl{Curl}

\DeclareMathOperator\Diff{Diff}
\DeclareMathOperator\Den{Den}

\newtheorem{theorem}{Theorem}

\newtheorem{proposition}[theorem]{Proposition}
\newtheorem*{problem*}{Problem}
{\theoremstyle{definition}

\newtheorem{remark}[theorem]{Remark}
\newtheorem*{remark*}{Remark}
}

\begin{document}

	\par\noindent {\LARGE\bf
	Selective decay for the rotating shallow-water equations with a structure-preserving discretization.
	\par}

{\vspace{4mm}\par\noindent {\bf R\"udiger Brecht$^\dag$, Werner Bauer$^\ddag$, Alexander Bihlo$^\dag$, Fran\c cois Gay-Balmaz$^\S$ and Scott MacLachlan$^\dag$ 
	} \par\vspace{2mm}\par}

{\vspace{2mm}\par\noindent {\it
		$^{\dag}$~Department of Mathematics and Statistics, Memorial University of Newfoundland,\\ St.\ John's (NL) A1C 5S7, Canada
}}
{\vspace{2mm}\par\noindent {\it
		$^{\ddag}$~Imperial College London, Department of Mathematics, 180 Queen’s Gate, London SW7 2AZ, United Kingdom.
}}
{\vspace{2mm}\par\noindent {\it
		$^{\S}$~\'Ecole Normale Sup\'erieure/CNRS, Laboratoire de M\'et\'eorologie Dynamique, Paris, France.
}}

{\vspace{2mm}\par\noindent {\it
		\textup{E-mail:} rbrecht@mun.ca, w.bauer@imperial.ac.uk, abihlo@mun.ca, gaybalma@lmd.ens.fr, smaclachlan@mun.ca 
	}\par}

\vspace{4mm}\par\noindent\hspace*{10mm}\parbox{140mm}{\small		
Numerical models of weather and climate critically depend on long-term stability of integrators for systems of hyperbolic conservation laws.  While such stability is often obtained from (physical or numerical) dissipation terms, physical fidelity of such simulations also depends on properly preserving conserved quantities, such as energy, of the system.  To address this apparent paradox, we develop a variational integrator for the shallow water equations that conserves energy, but dissipates potential enstrophy.  Our approach follows the continuous selective decay framework [F. Gay-Balmaz and D. Holm. Selective decay by Casimir dissipation in inviscid fluids. Nonlinearity, 26(2):495, 2013],  which enables dissipating an otherwise conserved quantity while conserving the total energy.  We use this in combination with the variational discretization method [D. Pavlov, P. Mullen, Y. Tong, E. Kanso, J. Marsden and M. Desbrun. Structure-preserving discretization of incompressible fluids. Physica D: Nonlinear Phenomena, 240(6):443-458, 2011] to obtain a \textit{discrete selective decay} framework. This is applied to the shallow water equations, both in the plane and on the sphere, to dissipate the potential enstrophy.  The resulting scheme significantly improves the quality of the approximate solutions, enabling long-term integrations to be carried out. 	 
}\par\vspace{4mm}

\section{Introduction}

Numerical weather and climate prediction requires the modeling of geophysical flows in the atmosphere and oceans on the globe. The atmosphere or ocean can be seen as a thin layer of fluid above the surface of the Earth and, thus, the shallow water equations (SWE) are a useful simplified model of the dynamics of geophysical flows around the Earth. These flows are approximately two-dimensional, and we can get insight into their flow dynamics by studying
the principles of two-dimensional turbulence. Important features of an
incompressible turbulent flow are the cascades of enstrophy and energy, where
the enstrophy transfers to higher wave numbers while the energy transfers to lower wave numbers, see e.g. Refs. \cite{lill71Ay,krai67Ay}. Numerical investigations of this phenomenon have led to the selective decay hypothesis \cite{matthaeus1980selective}, which states  that the enstrophy accumulation at the grid-scale should be dissipated to 
\textcolor{black}{damp this small scale noise} while energy should be conserved. 
\textcolor{black}{We note that, in this setting, energy conservation alone is  
	not enough to preclude the emergence of unphysical small scale oscillations (or: spurious modes) and to guarantee convergence to an accurate numerical solution.
	In practice, the dynamical core of a simulation needs to have several properties to guarantee accuracy, such as a control over spurious dynamical modes and mimetic and conservation properties, see Ref. \cite{Staniforth12} for details.}

Following \textcolor{black}{the selective decay hypothesis}, energy conserving and enstrophy dissipating numerical schemes for the SWE and other equations have been developed. 
Based on the anticipated vorticity method (APVM) \cite{sado85Ay}, an energy conserving and enstrophy dissipating finite-difference model was developed in Ref. \cite{arakawa1990energy}.
\textcolor{black}{Moreover, the anticipated vorticity model has been widely used to dissipate enstrophy and stabilize the fields, see Refs. \cite{mcra13Ay,li2020analysis,  ringler2010unified,ring08Ay}}.   
In Ref. \cite{shutts2005kinetic}, energy dissipated by a hyperviscosity term was ``reinserted'' to the system by adding random perturbations or specific large-scale velocity patterns.

\textcolor{black}{	
	In Ref. \cite{warneford2014thermal}, enstrophy was dissipated using a spectral filter.
	Further, in Ref. \cite{nair2021selective},  external forcing inputs are designed which change energy and enstrophy selectively.
}
Then, in Ref. \cite{thuburn2014cascades}, a scheme was proposed 
where the lost energy is repaired by adding a  vorticity perturbation to the  preliminary vorticity field at each timestep.
Moreover, in Ref. \cite{natale2017scale}, an energy-conserving and enstrophy-dissipating upwind-stabilization for finite-element discretizations was developed \textcolor{black}{which was extended and applied in Refs. \cite{wimmer2020energy, shipton2018higher}}.

However, many such methods directly manipulate the equations of motion to include the dissipation, which can have unpredictable consequences for the physical fidelity of the resulting numerical scheme. 
An alternative and more general method for developing energy-conserving dissipation schemes was introduced through the Lie--Poisson framework in Refs. \cite{gay2013selective,GBHo2014}.
While this approach appears to have great potential, it has not yet been applied to discrete models of geophysical fluid dynamics.  In this paper, we aim to ``connect the dots'', leveraging the energy conservation from the Lie--Poisson framework via a structure-preserving discretization method. Structure-preserving integrators for differential equations generally guarantee long-term stability, consistency in statistical properties, and  prevention of a systematic drift in stationary or periodic solutions, see Refs. \cite{hair06Ay,leim04Ay,wan16b}.

Here, we focus on variational integrators. These schemes are based on first discretizing the underlying variational principle and, then, deriving numerical schemes from the discrete Euler--Lagrange equations \cite{mars01a}.
In Ref. \cite{baue17a}, a variational discretization of the SWE was carried out. Following this,
in Ref. \cite{brecht2019variational}, the scheme was extended to the sphere, and it was observed that a stabilization of the scheme was needed to carry out long-term simulations of more than 50 days, \textcolor{black}{but also to avoid spurious small scale noise.}  
In this paper, we review the continuous selective decay theory to introduce a discretization of the selective decay that mimics the continuous theory. We apply the new framework to obtain a discretization of the SWE that dissipates enstrophy and conserves energy. In particular, we extend the discrete SWE introduced in Refs. \cite{baue17a, brecht2019variational} with the selective decay and carry out benchmarks in the plane and on the sphere.  

This article is structured as follows. In Section \ref{sec:background}, we review the continuous theory for variational discretization.  The continuous Casimir dissipation idea is introduced in Section \ref{sec:Selectivedecay}. Then, Section \ref{sec:discreteselective} is devoted to a description of the discrete integrator for the SWE proposed herein. In Section \ref{sec:numresults}, we verify the consistency of the discrete commutator and  present results from numerical simulations. The conclusions are given in Section \ref{sec:conclusion}. Furthermore, some detailed computations are presented in the appendix.

\section{Euler--Poincar\'{e} equations}\label{sec:background}

To obtain selective decay in the numerical scheme, we will use variational discretization, which mimics the continuous variational structure. On the continuous level, the equations of motion are obtained by defining a Lagrangian and computing the variational principle. This relies on the Euler--Poincar\'e reduction: the reformulation of Hamilton's principle from the Lagrangian to the Eulerian description. Thus, to understand the discretization procedure, we first review how we obtain the equations in the Euler--Poincar\'e framework.

The motion of a compressible fluid on a smooth manifold $ M $ (such as the surface of a sphere) is formally described by curves (functions) $\varphi :[0,T] \rightarrow \operatorname{Diff}( M )$ that are critical for the Hamilton principle,
\begin{equation}\label{eq:lagrangiantrajectory}
	\delta \int_0^T L(\varphi, \dot\varphi)~{\rm d} t=0,
\end{equation}
with respect to variations $\delta\varphi$ vanishing at $t=0$ and $t=T$. Here, $ \operatorname{Diff}( M )$ is the group of diffeomorphisms of the fluid domain $ M $ (differentiable one-to-one maps of $M$ onto itself with differentiable inverses), and $L$ is the Lagrangian of the fluid model expressed in terms of the Lagrangian fluid trajectory $\varphi$ and Lagrangian fluid velocity $\dot\varphi$. The variational principle \eqref{eq:lagrangiantrajectory} gives the equations in the Lagrangian description. For many computational approaches, it is more attractive to use a fixed Eulerian domain and, thus, Eulerian variables.
Rewriting the principle in Eq. \eqref{eq:lagrangiantrajectory} in Eulerian variables yields the Euler--Poincar\'e variational principle which involves constrained variations, see  Ref. \cite{holm98Ay} for a complete treatment.
Here, we give a brief overview and refer to the appendix of Ref. \cite{brecht2019variational} for a more detailed review for the case of the rotating shallow water equations on Riemannian manifolds.

We assume that $M$ is endowed with a Riemannian metric and denote by ${\rm d} \sigma $ the associated Riemannian volume form. The examples treated in this paper will be a doubly periodic domain in $ \mathbb{R} ^2 $ endowed with the Euclidean metric and a sphere endowed with its standard Riemannian metric; hence, we assume that $M$ has no boundary.
The Eulerian variables defined in terms of the Lagrangian fluid trajectory are the fluid velocity $\mathbf{u}=\dot\varphi\circ \varphi^{-1}\in \mathfrak{X}(M)$ (vector fields on $M$) and the fluid depth $h=\textcolor{black}{(h_0/J \varphi ) \circ \varphi ^{-1}  \in \Den(M)}$ (densities on $M$), where $h_0$ is the initial fluid depth and $J \varphi $ is the Jacobian of $ \varphi $ with respect to ${\rm d} \sigma $. The volume form allows the identification of the space of densities on $M$ with the space of functions on $M$.
From these relations,
the Lagrangian $L(\varphi , \dot \varphi )$ can be written in terms of $ \mathbf{u} $ and $h$, which yields the reduced Lagrangian $\ell\colon \mathfrak{X}(M)\times \Den(M)\to \R$. A consequence of the definition of $h$ is the mass continuity equation
\begin{equation}\label{mass_advection}
	\partial_t h + \operatorname{div}(h\mathbf{u})=0,
\end{equation}
with $\operatorname{div}$ being the divergence operator on $M$ defined by $ \pounds _ \mathbf{u} {\rm d} \sigma = (\operatorname{div} \mathbf{u} ) {\rm d} \sigma $.
Then, \eqref{eq:lagrangiantrajectory} yields the Euler--Poincar\'e variational principle with respect to constrained variations,
\begin{equation}\label{eq:varprincipal}
	\delta \int_0^T \ell(\mathbf{u},h) ~{\rm d} t=0 \quad \text{for}\quad 
	\begin{cases}
		~\delta \mathbf{u} = \partial_t \mathbf{v}+[\mathbf{u},\mathbf{v}]\\
		~\delta h = - \operatorname{div} (h\mathbf{v}),
	\end{cases}
\end{equation}
where $\mathbf{v}$ is an arbitrary \textcolor{black}{time-dependent} vector field with $\mathbf{v}(0)=\mathbf{v}(T)=0$ and $[ \cdot , \cdot ]$ is the Lie bracket of vector fields, $[\mathbf{u},\mathbf{v}]^i=\mathbf{u}^j \partial _j \mathbf{v}^i -\mathbf{v}^j \partial _j \mathbf{u}^i$ using the Einstein summation convention.

To compute the equations of motion in Eulerian variables, we need the functional derivatives $\frac{\delta \ell}{\delta \u}\in \Omega^1 (M)$ (one-forms on $ M $) and $\frac{\delta \ell}{\delta h}\in F(M)$  (scalar functions on $M$) which are defined by the duality pairings,
\begin{align}
	\left\langle \frac{\delta \ell}{\delta \u},\delta \u \right\rangle_1&:= \int_M \frac{\delta \ell}{\delta \u}\cdot \delta \u \, {\rm d} \sigma =\left. \frac{d}{d\varepsilon}\right|_{\varepsilon=0}  \ell(\mathbf{u}+\eps \delta \mathbf{u},h),
	\\
	\left\langle \frac{\delta \ell}{\delta h},\delta h \right\rangle_0&:= \int_M \frac{\delta \ell}{\delta h} \delta h \,  {\rm d} \sigma  = \left. \frac{d}{d\varepsilon}\right|_{\varepsilon=0}   \ell(\mathbf{u},h+\eps \delta h),
\end{align}
for arbitrary $ \delta \mathbf{u} $ and $ \delta h$.
Note that we denote the duality pairing between a one-form $\frac{\delta \ell}{\delta \u}$   and a vector field $\delta \u$  as $\langle \cdot,\cdot\rangle_1$ and the dual pairing between a function $\frac{\delta \ell}{\delta h}$ and a density $\delta h$ as $\langle \cdot,\cdot\rangle_0$.
The variational principle \eqref{eq:varprincipal} yields the Euler--Poincar\'e equations,
\begin{equation}\label{eq:EuPoin}
	\partial_t \frac{\delta \ell}{\delta \mathbf{u}} + \mathfrak{L} _{\mathbf{u}}\frac{\delta \ell}{\delta \mathbf{u}} = h \mathbf{d}\frac{\delta \ell}{\delta h},
\end{equation}
\textcolor{black}{where $ \mathfrak{L}  _{\mathbf{u}} \mathbf{m} = \mathbf{i} _ \mathbf{u} \mathbf{d} \mathbf{m}  + \mathbf{d} ( \mathbf{i} _ \mathbf{u} \mathbf{m} ) + \mathbf{m} \operatorname{div} \mathbf{u} $ for a one-form $ \mathbf{m} $ and a vector field $ \mathbf{u} $. Here $ \mathbf{i} _ \mathbf{u} \alpha $ denotes the contraction of a vector field $ \mathbf{u} $ with a differential form $ \alpha $ and $\d$ is the exterior derivative. We have the relation $ \pounds _ \mathbf{u} ( \mathbf{m} \otimes {\rm d} \sigma )=\mathfrak{L}_ \mathbf{u} \mathbf{m}  \otimes {\rm d} \sigma $, where $ \pounds _ \mathbf{u} $ is the Lie derivative of the one-form density $ \mathbf{m} \otimes {\rm d} \sigma $.}

\textcolor{black}{ In Euclidean space, using the identity
	$$
	\mathfrak{L}_{\mathbf{u}}\mathbf{v}= \mathbf{u}\cdot \nabla \mathbf{v} + \nabla \mathbf{u}^\top \mathbf{v}+\mathbf{v} \text{ div } \mathbf{u},
	$$
	the Euler--Poincar\'e equations reduce to 
	$$
	\partial_t \frac{\delta \ell}{\delta \mathbf{u}} + \mathbf{u}\cdot \nabla \frac{\delta \ell }{\delta \mathbf{u}} + \nabla\mathbf{u}^\top \frac{\delta \ell}{\delta \mathbf{u}}+\frac{\delta \ell}{\delta \mathbf{u}} \text{ div } \mathbf{u} = h \nabla \frac{\delta \ell}{\delta h}.
	$$}

\subsection{Variational principle for the SWE}
For the rotating shallow water equations on a two-dimensional Riemannian manifold $M$, the Lagrangian is given by
\begin{equation}\label{eq:lagrangSWEcont}
	\ell(\mathbf{u},h)=\int_M \Big[\frac{1}{2} h \mathbf{u}^\flat \cdot \mathbf{u}+h\mathbf{r}^\flat \cdot \mathbf{u}-\frac{1}{2}g(h+\eta_b)^2\Big]  {\rm d} \sigma ,
\end{equation}
where $\eta_b$ is the bottom topography, $g$ is the gravitational acceleration and $\mathbf{r}$ is the vector potential of the angular velocity of the Earth. Here, $\flat\colon TM\to TM^*$ is the flat operator of the Riemannian metric, that associates a one-form $\u^\flat$ to a vector field $\u$.
With the variational derivatives $\frac{\delta \ell}{\delta \mathbf{u}}=h(\mathbf{u}^\flat+\mathbf{r}^\flat)$ and $\frac{\delta \ell}{\delta h}=\frac{1}{2}\mathbf{u}^\flat \cdot \mathbf{u}+ \mathbf{r}^\flat\cdot\mathbf{u}-g(h+\eta_b)$, the Euler--Poincar\'e  equation \eqref{eq:EuPoin} gives the momentum equations of the SWE in the space of one-forms:
\begin{equation}\label{eq:RSWgeneral}
	\partial_t \mathbf{u}^\flat + \mathbf{i}_\mathbf{u} \mathbf{d} (\mathbf{u}^\flat + \mathbf{r}^\flat) + \mathbf{d}\left(\frac{1}{2}\mathbf{u}^\flat \cdot \mathbf{u} +g(h+\eta_b)\right)=0.
\end{equation}
This general expression reduces in the Euclidean space $\mathbb{R}^2$ to:
\begin{equation}\label{eq:RSW}
	\partial_t \mathbf{u} + \left(\nabla \times (\mathbf{u}+\mathbf{r})\right)\times \mathbf{u}+\nabla \left(\frac 1 2 |\mathbf{u}|^2+g(h+\eta_b)\right)=0.
\end{equation}

In the next section, we review a new dissipation scheme for our framework, that only acts on one conserved quantity while conserving the energy.  

\section{Selective decay with Casimir dissipation}\label{sec:Selectivedecay}

Given the Lagrangian $\ell(\u,h)$ of the fluid in Eulerian variables, the associated Hamiltonian function $\h(\m,h)$ is  obtained by the Legendre transformation
\[
\h (\m, h)=\langle \m, \u\rangle_1 -\ell(\u,h), 
\]
with $ \mathbf{u} $ defined in terms of $(\mathbf{m},h)$ by the relation $\m=\frac{\delta\ell}{\delta \u} \in \Omega ^1 (M)$
. We note the relations
\begin{equation}\label{eq:hamiltonRelations}
	\frac{\delta H}{\delta \m}=\u \quad\text{and}\quad \frac{\delta H}{\delta h}=-\frac{\delta \ell}{\delta h}.
\end{equation} 
The Eulerian Lie--Poisson formulation is given by
\begin{equation}\label{LP_form}
	\frac{df}{dt} =\{f,\h\}, \quad \forall\;f,
\end{equation}	
with Lie--Poisson bracket $\{ \cdot , \cdot \}$ defined as
\begin{align}\label{eq:liePoissonBracket}
		\{f,\h\}&=
		- \int_ M \mathbf{m} \cdot \Big[\frac{\delta f}{\delta \m},\frac{\delta \h}{\delta \m}\Big]  {\rm d} \sigma  
		 + \int_ M h  \Big( \mathbf{d} \frac{\delta f}{\delta h} \cdot \frac{\delta \h}{\delta \mathbf{m} } -   \mathbf{d} \frac{\delta \h}{\delta h} \cdot \frac{\delta f}{\delta \mathbf{m} }  \Big)  {\rm d} \sigma,
\end{align}
see Ref. \cite{holm98Ay} for details. The Lie--Poisson equations \eqref{LP_form} are equivalent to the system of equations \eqref{mass_advection} and \eqref{eq:EuPoin}, as can \textcolor{black}{be directly verified by} using \eqref{eq:hamiltonRelations}.

\paragraph{Example for the SWE.}
The Hamiltonian gives the total energy of the system; for the SWE, it reads
\begin{equation}
	H(\mathbf{m},h)=\int_M\Big[ \frac{1}{2h}|\mathbf{m}-h\mathbf{r}|^2+\frac{1}{2}g(h+\eta_b)^2\Big]{\rm d}\sigma.
\end{equation}
In this case, we have $\mathbf{m}=h(\mathbf{u}+\mathbf{r})$, so that $\frac{1}{2h}|\mathbf{m}-h\mathbf{r}|^2=\frac{1}{2}h|\mathbf{u}|^2$, which is the kinetic energy of the fluid.

For the selective decay, we use the relationship between the Lie--Poisson bracket and the conservation laws. A function $C$ is called a \textit{Casimir} for the Lie--Poisson bracket if it satisfies $\{C,f\}=0$ for all $f$.  With this, we have the conservation law $\frac{dC}{dt}=0$ along solutions of the Lie--Poisson system $\dot f= \{f, \h\}$, for any Hamiltonian $ \h$. In the next section, the Lie--Poisson bracket is extended to dissipate a Casimir but still conserve energy.  In Section~\ref{sec:contEnstCas}, we give a concrete example of a Casimir for the SWE.

\subsection{Casimir dissipation}

In this section, we recall the approach to selective decay developed in Ref. \cite{gay2013selective}.
Let $\gamma\colon \mathfrak{X}(M)\times \mathfrak{X}(M)\to \R$ be a positive and symmetric bilinear form (with associated norm $\|\cdot\|_\gamma$) and $C$ a Casimir function. 
The Casimir dissipation is introduced in the Lie--Poisson formulation as follows
\begin{equation}\label{eq:casimirDissipation}
	\frac{df}{dt}=\{ f,\h \}-\theta \gamma\left(\Big[\frac{\delta f}{\delta \m},\frac{\delta \h}{\delta \m}\Big], \Big[ \frac{\delta C}{\delta \m},\frac{\delta \h}{\delta \m}\Big]\right),
\end{equation}
for some $ \theta >0$.
If $f=\h$, then (from the definition of the Lie bracket) $\left[\frac{\delta \h}{\delta \m},\frac{\delta \h}{\delta \m}\right]=0$, giving
$$
\frac{d\h}{dt}=0-\theta\gamma\left(0,\Big[\frac{\delta C}{\delta \m},\frac{\delta \h}{\delta \m}\Big]\right)=0,
$$
and we see that the energy remains conserved. For $f=C$, we have $$
\frac{dC}{dt}=0-\theta\gamma\left(\Big[\frac{\delta C}{\delta \m},\frac{\delta \h}{\delta \m}\Big],\Big[\frac{\delta C}{\delta \m},\frac{\delta \h}{\delta \m}\Big]\right)=-\theta \left\|\Big[\frac{\delta C}{\delta \m},\frac{\delta \h}{\delta \m}\Big]\right\|^2_\gamma,
$$
thus, the Casimir decays in time when $\left[\frac{\delta C}{\delta \m},\frac{\delta \h}{\delta \m}\right] \not\equiv 0$.

The corresponding Lagrange--d'Alembert variational principle is given by (see Ref. \cite[ Eq. (3.7)]{gay2013selective})
\begin{align}\label{eq:varprincCas}
	\begin{split}
		\Big(\delta\int_0^T\ell(\mathbf{u},h) {\rm d} t \Big) &+\theta \int_0^T \gamma\left( \Big[\frac{\delta C}{\delta \mathbf{m}},\mathbf{u}\Big],\Big[\mathbf{u},\mathbf{v}\Big]\right) {\rm d} t=0, 
		\\
		&\quad\text{for}\quad 
		\begin{cases}
			~\delta \mathbf{u} = \partial_t \mathbf{v}+[\mathbf{u},\mathbf{v}]\\
			~\delta h = - \operatorname{div} (h\mathbf{v}).
		\end{cases}
	\end{split}
\end{align}
Then, the Casimir dissipative Euler--Poincar\'e equations (see Ref. \cite[ Eq. (3.3)]{gay2013selective}) are
\begin{equation}\label{eq:EuPoinCas}
	\partial_t \frac{\delta \ell}{\delta \mathbf{u}} + \textcolor{black}{\mathfrak{L} _{\mathbf{u}}}\frac{\delta \ell}{\delta \mathbf{u}} = h \mathbf{d}\frac{\delta \ell}{\delta h}+\theta \textcolor{black}{\mathfrak{L}_ \mathbf{u}} \Big( \Big[\mathbf{u},\frac{\delta C}{\delta \mathbf{m}}\Big]^ \gamma  \Big),
\end{equation}
where, for a vector field $ \mathbf{u} \in \mathfrak{X} ( M )$, $ \mathbf{u} ^ \gamma $ is the one-form on $ M $ defined by $\int_ M (\mathbf{u} ^\gamma  \cdot \mathbf{v}) {\rm d} \sigma = \gamma ( \mathbf{u} , \mathbf{v} )$, for all $\mathbf{v} \in \mathfrak{X} ( M )$, and we recall that $ \mathbf{m} = \frac{\delta \ell}{\delta \mathbf{u} }$. We assume that $ \gamma $ is such that the one-form $ \mathbf{u} ^ \gamma $ is well-defined for all $ \mathbf{u} \in \mathfrak{X} (M)$, see Ref. \cite{gay2013selective} for examples.

\subsection{Enstrophy dissipation for SWE}\label{sec:contEnstCas}
Next, we will consider enstrophy dissipation for the SWE.
For two-dimensional fluid flows dominated by geostrophic balance, enstrophy is known to cascade to small scales.  Thus,
in order to obtain physically relevant solutions, it is necessary to dissipate enstrophy at such scales, see Refs. \cite{bona05Ay,mcra13Ay,ring02a}. For the SWE on two-dimensional Riemannian manifolds, the potential enstrophy Casimir is given by
\begin{equation}
	C(\mathbf{m},h)= \frac 1 2 \int_M h \,q(\mathbf{m},h)^2 {\rm d} \sigma
\end{equation} 
\textcolor{black}{with $q( \mathbf{m} , h)$ the potential vorticity function defined by}
\[ 
q(\mathbf{m},h) {\rm d} \sigma = \frac{1}{h} \mathbf{d} \frac{ \mathbf{m} }{h},
\]
\textcolor{black}{where we recall that the 2-form $ {\rm d} \sigma $ is
	the Riemannian volume form.}
The variational derivative of the enstrophy Casimir is found as
$\frac{\delta C}{\delta \m} = -\frac{1}{h} (\star \mathbf{d}  q)^\sharp$ with $\sharp:T^*M \rightarrow TM$ the Riemannian sharp operator, see Appendix \ref{sec:VarDerEnstCas}. For a two-dimensional planar domain, these formulas reduce to
\begin{equation} 
	q(\mathbf{m},h)= \frac{1}{h} \mathbf{z}\cdot \nabla\times\left(\frac{\mathbf{m}}{h} \right), \qquad \frac{\delta C}{\delta \m} = -\frac{1}{h} \mathbf{z}\times\nabla q,
\end{equation}
where $\mathbf{z}$ is the canonical unit vector pointing in the positive $z$-direction.
With the Lagrangian \eqref{eq:lagrangSWEcont}, the Casimir dissipative Euler--Poincar\'e equations \eqref{eq:EuPoinCas} are given by 
\begin{align}\label{eq:RSWcasdis}
	\begin{split}
		h\partial_t \mathbf{u}^\flat + h\mathbf{i}_\mathbf{u} \mathbf{d} (\mathbf{u}^\flat + \mathbf{r}^\flat)
		&=-h \mathbf{d}\left(\frac{1}{2}\mathbf{u}^\flat \cdot \mathbf{u} +g(h+\eta_b)\right)
		\\
		&~~~~ +\theta \textcolor{black}{\mathfrak{L} _{\mathbf{u}}}\Big(h\Big[\mathbf{u},\frac{\delta C}{\delta \mathbf{m}}\Big]^\flat\Big),	
	\end{split}
\end{align}
where we choose $\gamma$ to be the water depth weighted $L^2$ inner product, i.e., $ \gamma ( \mathbf{u} , \mathbf{v} ) =\int_ M h (\mathbf{u} ^\flat \cdot \mathbf{v})  {\rm d} \sigma $, and we note that $ \mathbf{u} ^ \gamma = h \mathbf{u} ^\flat $, with $ \flat$ associated to the Riemannian metric on $ M $.
In $\R^2$ the Casimir dissipating equation reduces to
\begin{align}\label{eq:RSWenst}
	\begin{split}
		\partial_t \mathbf{u} + &\left(\nabla \times (\mathbf{u}+\mathbf{r})\right)\times \mathbf{u}+\nabla \left(\frac 1 2 |\mathbf{u}|^2+g(h+\eta_b)\right)
		\\
		&\qquad \qquad =\theta\left( \mathbf{u}\cdot \nabla \mathbf{w} + \nabla \mathbf{u}^\top \mathbf{w}+\mathbf{w} \text{ div } \mathbf{u}\right),
	\end{split}
\end{align}
where $\mathbf{w}=[\mathbf{u},-\frac{1}{h} \mathbf{z}\times\nabla q]$.

\section{Discrete selective decay}\label{sec:discreteselective}

The discretization process translates each step of the continuous theory to the discrete level.  Here, we review the variational discretization process for fluid initially developed in Ref. \cite{pavlov2011structure}, see also Refs. \cite{baue17a,desbrun2014variational,GaGB2019,gawlik2011geometric} for extensions, and incorporate the Casimir  selective decay into it.

\begin{figure}
	\centering
	\includegraphics[width=0.5\textwidth]{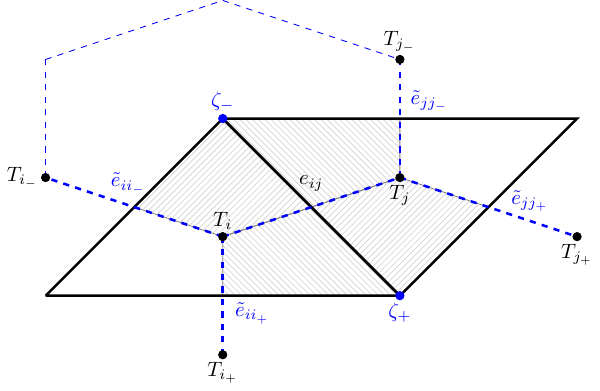}
	\caption{Notation and indexing conventions for the two-dimensional simplicial mesh.}
	\label{fig:Notationgrid}
\end{figure}

We consider a two-dimensional simplicial mesh $\M$  with $n$ cells on the fluids domain, where  triangles ($T$) are used as the primal grid, and the circumcenter dual ($\zeta$) as the dual grid. On the grid (see Fig. \ref{fig:Notationgrid}) we adopt the following notation:
\begin{itemize}
	\item $e_{ij}=T_i\cap T_j$ as the primal edge,
	\item $\tilde{e}_{ij}= \zeta_{+}\cap \zeta_{-}$ as the dual edge,
	\item $\Omega_{ii}$ as the area of triangle $T_i$,
	\item $h_i$ as the discrete water depth on $T_i$,
	\item $(\eta_b)_i$ as the discrete bottom topography on $T_i$,
	\item $V_{ij}$ is $\mathbf{u}_{e_{ij}}\cdot\mathbf{n}_{e_{ij}}$ at the edge midpoint. 
	\item $\overline{h}_{ij}=\frac{1}{2}(h_i+h_j)$	as the water depth averaged to the edge midpoints.
\end{itemize}
Here, $\mathbf{n}_{e_{ij}}$ is the normal vector of edge $e_{ij}$ pointing towards $T_j$.

\subsection{Discrete setup}

The Euler--Poincar\'e reduction %
from the Lie group of diffeomorphisms to its Lie algebra (as discussed in Section \ref{sec:background}) is done analogously in the discrete setting by identifying the discrete analogues of $\Diff(M)$ and $\Den(M)$.
For piecewise constant functions, the discrete diffeomorphism group is the matrix group
\begin{equation}
	D(\M)=\{ q\in GL(n)^+ \mid q\cdot \mathbf{1}=\mathbf{1} \},
\end{equation}
with $GL(n)^+$ the group of real $n \times n$ matrices with positive determinant and $\mathbf{1}=(1,...,1)^\mathsf{T}$.
The condition $q\cdot \mathbf{1}=\mathbf{1}$ means that constants are preserved, which is needed to obtain mass conservation. Then, the Lie algebra of $D(\M)$ is
\begin{equation}
	\mathfrak{d}(\M)=\{ A\in \mathfrak{gl}(n) \mid A\cdot \mathbf{1}=0\}
\end{equation}
with the matrix commutator $[A,B]=AB-BA$ as the Lie bracket, where $\mathfrak{gl}(n)$ is the Lie algebra of $n\times n$ real matrices.
\textcolor{black}{To describe the infinitesimal exchanges of fluid particles between adjacent pairs of cells, a subspace $ \mathcal{R} \subset \mathfrak{d}(\M)$ is considered \cite{baue17a}, which corresponds to discrete vector fields. This subspace is given by
	\begin{align*}
		R_1 &= \big\{ A \in \mathfrak{d}(\M)\mid A^\top \Omega + \Omega A^\top \text{ is diagonal }\big\}
		\\
		R_2&=\big\{ A \in \mathfrak{d}(\M)\mid A_{ij}=0 ~~~\forall j\notin N(i)\big\}
		\\
		\mathcal{R}&= R_1 \cap R_2
		,
	\end{align*}
	with $N(i)$ being the set of cells sharing an edge with the cell $T_i$ and with $ \Omega$ being the $n \times n$ diagonal matrix with elements $ \Omega _{ii}$.}

\begin{remark}\label{r:spaceR}
	For $A,B\in \mathcal{R}$ we have $[A,B]_{ij}=0$ for all $j \in N(i)$. Since elements of $\mathcal{R}$ are zero for non-neighboring cells, we get $[\mathcal{R},\mathcal{R}]\cap \mathcal{R}=\{0\}$. In particular $[ \mathcal{R} , \mathcal{R} ] \neq \mathcal{R} $ hence the subspace $ \mathcal{R}  \subset \mathfrak{d}(\M)$ corresponds to a nonholonomic constraint. 
Consequently, we need to define a discrete commutator $\mybrL , \mybrR$ such that $\mybrL A, B\mybrR \in \mathcal{R}$, so that we can directly apply the definition of the discrete flat operator to the discrete commutator. In section \ref{sec:casdissscheme}, we give the details on how to obtain such a discrete commutator. 
\end{remark}

Next, we identify the dual space $\mathcal{R}^*$ with the space $\Omega^1_d(\M)$ of discrete one-forms relative to the duality pairing on $\mathfrak{gl}(n)$:
\begin{equation}\label{eq:duality1}
	\langle L, A\rangle_1 = \text{Tr}(L^\top \Omega A).
\end{equation}
To obtain an element in $\mathcal{R}^*$, we use the projection $P\colon \mathfrak{gl}(n)\to \Omega^1_d(\M)$ defined by
\begin{equation}\label{eq:Discreteprojection}	
	P(L)_{ij}=\frac{1}{2}(L_{ij}-L_{ji}-L_{ii}+L_{jj}),
\end{equation}
which satisfies $ \left\langle L, A \right\rangle _1 = \left\langle P(L), A \right\rangle _1$, for all $A \in \mathcal{R}$, see Ref. \cite{baue17a}.
Piecewise constant functions on $M$ are represented by vectors $F\in\R^n$, with value $F_i$ on cell $i$ being the cell average of the continuous function on cell $i$. The space of discrete functions is denoted by $\Omega^0_d(\M)$, and the space of discrete densities $\Den_d(\M)\simeq \R^n$ is defined as the dual space to $\Omega^0_d(\M)$ relative to the pairing:
\begin{equation}
	\langle F, G\rangle_0 = F^\top \Omega G.
\end{equation}
If a matrix $A\in \mathcal{R}$ approximates a vector field $\u$, then its entries satisfy
\begin{align}\label{eq:approxVel}
	\begin{split}
		A_{ij}&\approx -\frac{1}{2\Omega_{ii}}\int_{T_i\cap T_j} (\u\cdot \mathbf{n}) {\rm d} S, \;\; j\in N(i)
		\\
		A_{ii}&\approx \frac{1}{2\Omega_{ii}}\int_{T_i}  \operatorname{div}  \u \,{\rm d} \sigma .
	\end{split}
\end{align}
The discrete flat operator on $ \mathcal{R} $ is defined by the two conditions
\begin{align}\label{eq:discrFlat}
	\begin{split}
		A_{ij}^\flat &= 2\Omega_{ii}\frac{|\tilde{e}_{ij}|}{|e_{ij}|} A_{ij}, \qquad j\in N(i)\\
		A_{ij}^\flat+A_{jk}^\flat+A_{ki}^\flat &= K^\zeta_j \langle \omega(A^\flat), \zeta\rangle, \qquad i, k\in N(j), k\notin N(i),
	\end{split}
\end{align}
where the cells $i, j, k$ have a common node, whose dual cell is $\zeta$, where
$K^\zeta_j=\frac{|\zeta\cap T_j|}{|\zeta|}$ and
$$
\langle\omega(A^\flat), \zeta\rangle=\sum_{\tilde{e}_{nm}\in\partial \zeta} A^\flat_{nm}
$$
is the discrete vorticity at cell $\zeta$. The definition of $A_{ij}^\flat$ represents the flux $2\Omega_{ii} A_{ij}$ multiplied by the Hodge star. The definition for non-adjacent cells uses that the exterior derivative \textcolor{black}{ $(\mathbf{d} A^\flat)_{ijk}$ is} a fraction of the vorticity at the dual node. \textcolor{black}{For details about the derivation of the discrete flat operator and the resulting weights, we refer the reader to Ref. \cite{pavlov2011structure}}.

\textcolor{black}{
	\begin{remark}
		As shown in Ref. \cite{baue17a}, standard assumptions on the mesh are needed to show that the approximations $A\in \mathcal{R}$ converge to continuous vector fields. \textcolor{black}{In particular, we assume that the mesh belongs to a shape-regular, quasi-uniform family of triangulations of $M$.}
	\end{remark}
}

In the next section, we will use this discrete setup to state the discrete variational principle and compute the numerical scheme.

\subsection{Discrete variational equations for selective decay}
Let $\ell\colon  \mathfrak{d}(\M) \times \Den_d(\M)\to \R$ be a semi-discrete Lagrangian and $C\colon  \mathfrak{d}(\M)  \times \Den_d(\M)\to \R$ be a semi-discretized approximation of a Casimir.
As above, let  $\gamma\colon  \mathfrak{d}(\M)\times  \mathfrak{d}(\M)\to \R$ be a positive, symmetric bilinear form. 
Analogous to the continuous Casimir dissipative variational principle in Eq. \eqref{eq:varprincCas}, we consider the discrete dissipative variational principle given by
\begin{align}\label{eq:discvarprinCas}
	\begin{split}
		\Big(\delta \int_0^T \ell(A,h) {\rm d} t \Big)&+\theta\int_0^T \gamma\left( \Big{\mybrL}\frac{\delta C}{\delta M},A\Big{\mybrR},[A,B] \right){\rm d} t=0\qquad 
		\\
		&\text{for} \qquad
		\begin{cases}
			\delta A = \partial_t B+[B, A]\\
			\delta h = -\Omega ^{-1}  B^\top \Omega h
		\end{cases}
	\end{split}
\end{align}
\textcolor{black}{where $A(t) \in \mathcal{R} $ and $B(t)$ is an arbitrary curve in $\mathcal{R}$ with $B(0)=B(T)=0$. This means that the constraint $\mathcal{R} \subset \mathfrak{d}(\M)$ is treated as a nonholonomic constraint in the variational principle, exactly as in Refs. \cite{pavlov2011structure,baue17a}.}
The discrete functional derivatives $\frac{\delta \ell}{\delta A}\in \mathfrak{d}(\M) ^*, \frac{\delta \ell}{\delta h}\in \Omega^0(\M)$ and $\frac{\delta C}{\delta M}\in \mathfrak{d}(\M) $ are defined by
\begin{align}\label{eq:varderivative}
	\begin{split}		
		\left\langle \frac{\delta \ell}{\delta A}, \delta A \right\rangle_1 &= \left. \frac{d}{d\varepsilon}\right|_{\varepsilon=0} \ell(A+\eps \delta A, h),
		\\
		\left\langle \frac{\delta \ell}{\delta h}, \delta h \right\rangle_0 &= \left. \frac{d}{d\varepsilon}\right|_{\varepsilon=0} \ell(A, h+\eps \delta h),
		\\
		\left\langle  \delta M, \frac{\delta C}{\delta M} \right\rangle_1 &= 
		\left. \frac{d}{d\varepsilon}\right|_{\varepsilon=0}  C(M+\eps \delta M, h),
	\end{split}
\end{align}
for all $ \delta A \in \mathfrak{d}(\M)  $, $ \delta h \in \operatorname{Den}_d( \mathbb{M})$, $ \delta M \in  \mathfrak{d}(\M) ^* $.

\begin{theorem}[\textbf{Discrete dissipative variational equations}]\label{th:discreteVarEq}
	For a semi-discrete Lagrangian $\ell(A,D)$, the curves $A(t),h(t) \textcolor{black}{\in \mathcal{R}}$ are critical for the variational principle of Eq. \eqref{eq:discvarprinCas} if and only if they satisfy
	\begin{equation}	\label{eq:projectionCas}
		P \left(	\frac{d}{dt}\frac{\delta \ell}{\delta A}+ \mathcal{L}_A\Big(\frac{\delta \ell}{\delta A}\Big)- \theta\mathcal{L}_A\Big(h  \Big{\mybrL}\frac{\delta C}{\delta M},A\Big{\mybrR}^\flat\Big)+ h\frac{\delta \ell}{\delta h}^\top\right)_{ij}=0,
	\end{equation}				
	\textcolor{black}{where $\mathcal{L}$ is the discrete analog to $\mathfrak{L}$ and it is defined by the commutator via the following relation 
		\begin{equation} \label{eq:liederdef}
			\left\langle \mathcal{L}  _A M, B \right\rangle _1 = \left\langle M, [A,B] \right\rangle _1.
		\end{equation}
	}	
\end{theorem}
\textbf{Proof:} 
The variational principle \eqref{eq:discvarprinCas} gives
\begin{align*}
	0&=
	\delta \int_0^T \ell(A,h){\rm d} t+\theta\int_0^T \gamma\left( \Big{\mybrL}\frac{\delta C}{\delta M},A\Big{\mybrR} ,[A,B] \right){\rm d} t.
	\intertext{Next, we use the definition of the flat operator and $\gamma$ to be the water depth weighted inner product, giving}
	0&=\delta \int_0^T \ell(A,h){\rm d} t+\theta\int_0^T \left\langle h \Big{\mybrL}\frac{\delta C}{\delta M},A\Big{\mybrR}^\flat,[A,B] \right\rangle_1{\rm d} t.
	\intertext{Finally, we use the expression of the variations in \eqref{eq:discvarprinCas} and \textcolor{black}{the property of Eq. \eqref{eq:liederdef},} which yields}
	0 &=
	-\int_0^T \left\langle
	\frac{d}{dt}\frac{\delta \ell}{\delta A}+ \mathcal{L}_A\Big(\frac{\delta \ell}{\delta A}\Big)+h\frac{\delta \ell}{\delta h}^\top 
	,B \right \rangle_1 {\rm d} t
	\\
	&~~~~ +\theta\int_0^T \left\langle h \Big{\mybrL}\frac{\delta C}{\delta M},A\Big{\mybrR}^\flat,[A,B] \right\rangle_1{\rm d} t
	\\&=
	-\int_0^T \left\langle
	\frac{d}{dt}\frac{\delta \ell}{\delta A}+ \mathcal{L}_A\Big(\frac{\delta \ell}{\delta A}\Big)+h\frac{\delta \ell}{\delta h}^\top 
	,B \right \rangle_1 {\rm d} t
	\\
	&~~~~ +\theta\int_0^T \left\langle \mathcal{L}_A\Big( h \Big{\mybrL}\frac{\delta C}{\delta M},A\Big{\mybrR}^\flat\Big),B \right\rangle_1	{\rm d} t
	\\&=
	-\int_0^T \Big\langle
	\frac{d}{dt}\frac{\delta \ell}{\delta A}+ \mathcal{L}_A\Big(\frac{\delta \ell}{\delta A}\Big)
	\\
	&~~~~ -\theta \mathcal{L}_A\Big(  h \Big{\mybrL}\frac{\delta C}{\delta M},A\Big{\mybrR}^\flat\Big)+h\frac{\delta \ell}{\delta h}^\top 
	,B \Big \rangle_1{\rm d} t.
\end{align*}
The result then follows from $\int_0^T \langle L, B\rangle_1~{\rm d} t=0,~ \forall B\in \mathcal{R} \iff  P(L)_{ij}=0$ (see Ref. \cite[Proposition 2.3]{baue17a}). 

\textcolor{black}{
	Note that the calculation above assumes that $\Big{\mybrL}\frac{\delta C}{\delta M},A\Big{\mybrR}$ is represented as a discrete vector field in $\mathcal{R}$, so that we can use the definition of the discrete flat operator given in \eqref{eq:discrFlat}.  As we will see in Section \ref{sec:casdissscheme}, this is \textcolor{black}{made possible by replacing the} direct calculation of the vector field commutator for elements of $\mathcal{R}$ with a representation using vector calculus identities.
}
\begin{flushright}
	$\square$	
\end{flushright}

\begin{remark}  We note that Theorem \ref{th:discreteVarEq} becomes the \emph{discrete variational equations} theorem  in Ref. \cite{baue17a} for $\theta=0$. This form of the discrete equations is valid on Cartesian and simplicial meshes in two and three dimensions. We focus below on two-dimensional simplicial meshes.
\end{remark}
The following proposition demonstrates that, for the resulting semi-discrete scheme, the energy is conserved.
\begin{proposition}\label{prop:compEnergy}
	Let $A(t)$ and $h(t)$ be the solution of \eqref{eq:projectionCas} and $ \dot h + \Omega ^{-1} A^\top \Omega h=0$.  Then,
	$$
	\frac{d}{dt}\left( \left\langle  \frac{\delta \ell}{\delta A},A \right\rangle _1 -\ell(A,h) \right)=0
	$$
\end{proposition}
\textbf{Proof:} 
We compute
\begin{align*}
	&\frac{d}{dt}\left( \left\langle  \frac{\delta \ell}{\delta A},A \right\rangle _1 -\ell(A,h) \right)
	\\
	&= \left\langle \frac{d}{dt}  \frac{\delta \ell}{\delta A},A \right\rangle _1 + \left\langle  \frac{\delta \ell}{\delta A}, \frac{d}{dt} A \right\rangle _1- \left\langle  \frac{\delta \ell}{\delta A}, \frac{d}{dt} A \right\rangle _1 - \left\langle  \frac{\delta \ell}{\delta h}, \frac{d}{dt} h \right\rangle _0
	\\
	&=  \left\langle P\frac{d}{dt}  \frac{\delta \ell}{\delta A},A \right\rangle _1 + \left\langle  \frac{\delta \ell}{\delta h},  -\Omega ^{-1} A^\top \Omega h \right\rangle _0
	\\
	&=  \left\langle P \left( \frac{d}{dt}  \frac{\delta \ell}{\delta A} + h \frac{\delta \ell}{\delta h}^\top \right) , A \right\rangle _1
	\\
	& = - \left\langle P \left( \mathcal{L} _ A \frac{\delta \ell}{\delta A} - \theta \mathcal{L} _A \Big(h  \Big{\mybrL}\frac{\delta C}{\delta M},A\Big{\mybrR}^\flat\Big)\right) , A \right\rangle_1 
	\\
	& = - \left\langle   \mathcal{L} _ A \frac{\delta \ell}{\delta A} - \theta \mathcal{L} _A \Big(h  \Big{\mybrL}\frac{\delta C}{\delta M},A\Big{\mybrR}^\flat\Big) , A \right\rangle _1=0 ,
\end{align*} 
where the last equality follows from \textcolor{black}{the \textcolor{black}{definition} in Eq. \eqref{eq:liederdef}.} %
This holds independently of the chosen discretization of $\Big{\mybrL}\frac{\delta C}{\delta M},A\Big{\mybrR}^\flat$.
\begin{flushright}
	$\square$	
\end{flushright}

\subsection{Variational discretization of the Casimir dissipative SWE}

Before presenting the discretization including the Casimir dissipation term, we briefly recall the variational discretization for the scheme without Casimir dissipation \cite{baue17a,brecht2019variational}.
We discretize the Lagrangian  \eqref{eq:lagrangSWEcont} with piecewise constant functions, giving
\begin{align}\label{eq:discrteLagrangian}
	\begin{split}
		\ell(A, h) &= \frac{1}{2}\sum_{i,j=1}^n h_i A^\flat_{ij}A_{ij} \Omega_{ii} 
		+\sum_{i,j=1}^n h_i R^\flat_{ij} A_{ij} \Omega_{ii}
		\\
		&~~~~ -\frac{1}{2}\sum_{i=1}^n g(h_i+(\eta_b)_i)^2\Omega_{ii}.
	\end{split}
\end{align} 
To compute the variational derivatives, we use the duality pairing \eqref{eq:duality1} and the definition in Eq. \eqref{eq:varderivative}, giving
\begin{align}\label{eq:funcder}
	\begin{split}
		\frac{\delta \ell}{\delta A}_{ij}&=h_i (A_{ij}^\flat+R_{ij}^\flat)
		\\
		\frac{\delta \ell}{\delta h}_i&=\frac{1}{2}\sum_j A^\flat_{ij}A_{ij}+\sum_j R^\flat_{ij}A_{ij}-g(h_i+(\eta_b)_i).
	\end{split}
\end{align}

In Refs. \cite{brecht2019variational,brecht2021rotating}, it was noted that the approximations of the differential operators that result from the variational discretization method agree with the following standard finite difference and finite volume operators:
\begin{align}\label{eq:discOps}
	\begin{split}
		(\text{Grad}_n~F)_{ij}&:=\frac{F_{T_j}-F_{T_i}}{|\tilde{e}_{ij}|},
		\\
		(\text{Grad}_t~F)_{ij}&:=\frac{F_{\zeta_-}-F_{\zeta_+}}{|e_{ij}|},
		\\
		(\operatorname{div} \mathbf{u})_i\approx		(\text{Div}~V)_i&:=\frac{1}{\Omega_{ii}}\sum_{k\in\{j,i_-,i_+\}}|e_{ik}|V_{ik},
		\\
		(\nabla \times \mathbf{u})_\zeta\approx	(\text{Curl}~V)_{\zeta}&:=\frac{1}{|\zeta|}\sum_{\tilde{e}_{nm}\in\partial\zeta} |\tilde{e}_{nm}|V_{nm} ,
	\end{split}
\end{align}
for a scalar field $F$ sampled either at the triangle or dual cell centres and a normal velocity $V_{ij}$. The components of the gradient of $F$ in the tangential and normal directions to an edge are denoted by $\text{Grad}_t$ and $\text{Grad}_n$, respectively. The normal velocity $V_{ij}$ is related to the matrix elements $A \in \mathcal{R}$ in \eqref{eq:approxVel} as 
\[
A_{ij}= -\frac{|e_{ij}|}{2\Omega_{ii}}  V_{ij}, \; j\in N(i) \quad\text{and}\quad A_{ii} =  \frac{1}{2\Omega_{ii}}\sum_{k\in N(i)} |e_{ik}|V_{ik}.
\]
\begin{remark}\label{r:curlmap} The gradient in the normal direction and the discrete divergence are adjoints with respect to the natural inner products on the triangles and their edges.  Similarly, the tangential gradient and the discrete curl operator are adjoints with respect to the natural inner products on dual cells and their edges.
\end{remark}
\begin{remark}
	The continuous gradient, divergence, and curl operators are naturally written in Cartesian coordinates, but can also be defined (via parametrization) in a local neighbourhood on the sphere. The discrete counterparts are always locally defined and independent of the coordinate system. We will use the notation in Cartesian  coordinates for this section for the continuum operators, to simplify notation.
\end{remark}
Computing the projection \eqref{eq:projectionCas} with $\theta=0$, we obtain the momentum equation in Ref. \cite{baue17a}.
For simplicity of presentation, we group the terms involved in the advection term $(\nabla\times (\mathbf{u}+\mathbf{r}))\times \mathbf{u}$ and denote their discretization by Adv$(V,h)$.  Similarly, the terms involved in the kinetic energy term $\nabla(\frac{1}{2}\mathbf{u}^2)$ are denoted by K$(V)$, and the terms involved in the gradient term $g\nabla h$ by G$(h)$.  Thus, we write
\begin{equation}\label{eq:dicsreteMomentumcas}
	\partial_t V_{ij}=-\text{Adv}(V,h)_{ij}-\text{K}(V)_{ij}-\text{G}(h)_{ij},
\end{equation}
where
\begin{align*}
	& \operatorname{Adv}_{\rm }(V,h)_{ij} := \\ 
	&   -  \frac{1}{\overline{h}_{ij} |\tilde{e}_{ij}| } \Big((\text{Curl } V)_{\zeta_{-}}+f_{\zeta_{-}}\Big) 
	\Big(   \frac{|\zeta_{-} \cap T_i|}{2 \Omega_{ii}  } \overline{h}_{ji_-} |e_{ii_-}|   V_{ii_-} \\
	& \hspace{12.5em} +  \frac{|\zeta_{-}  \cap T_j|}{2 \Omega_{jj}  } \overline{h}_{ij_-}  |e_{jj_-}|   V_{jj_-} \Big) \\
	& +  \frac{1}{\overline{h}_{ij}  |\tilde{e}_{ij}| }\Big((\text{Curl } V)_{\zeta_{+}}+f_{\zeta_{+}}\Big) 
	\Big(  \frac{|\zeta_{+}  \cap T_i|}{2 \Omega_{ii} }\overline{h}_{ji_+} |e_{ii_+}|   V_{ii_+} \\
	& \hspace{12.5em} + \frac{|\zeta_{+}  \cap T_j|}{2 \Omega_{jj}   } \overline{h}_{ij_+}   |e_{jj_+}|   V_{jj_+} \Big),
	\\
	& \operatorname{K}_{\rm }( V) _{ij} := \frac{1}{2}(\text{Grad}_n ~F)_{ij}, \\
	& F_{T_i}= \sum_{k\in \{j, i_-, i_+\} }\frac{|\tilde{e}_{ik}|~|e_{ik}|(V_{ik})^2}{2\Omega_{kk}},
	\\
	&\operatorname{G}(h)_{ij} :=  g (\text{Grad}_n~(h+\eta_b))_{ij}.
\end{align*}
The Coriolis parameter is defined by
$$
f_\zeta=\frac{1}{|\zeta|}\sum_{\tilde{e}_{nm}\in \partial \zeta}|\tilde{e}_{nm}|r_{nm},
\qquad \text{ with } r_{ij}=\mathbf{r}_{e_{ij}}\cdot \mathbf{n}_{e_{ij}},
$$
\textcolor{black}{where $\mathbf{r}_{e_{ij}}$ is the vector potential of the angular velocity of the Earth evaluated at the edge midpoint.}

\subsection{Casimir dissipative scheme}\label{sec:casdissscheme}

Including the extra term for $\theta>0$, the Casimir dissipative momentum equation is
\begin{equation}\label{eq:momentumCas}
	\partial_t V_{ij}=-\text{Adv}(V,h)_{ij}-\text{K}(V)_{ij}-\text{G}(h)_{ij} +\theta \text{L}(V,h,\frac{\delta C}{\delta M})_{ij},
\end{equation}
where
$$
\frac{1}{h_i|\tilde{e}_{ij}|}P\left(\mathcal{L}_A\Big(h  \Big{\mybrL}\frac{\delta C}{\delta M},A\Big{\mybrR}^\flat\Big)\right)_{ij}=:\text{L}(V,h,\frac{\delta C}{\delta M})_{ij}
$$ 
\textcolor{black}{that follows from a comparison between \eqref{eq:projectionCas} to \eqref{eq:momentumCas}.}

To compute the latter term, we first need to discretize the commutator.
Here, we cannot follow the discretization procedure of Ref. \cite{pavlov2011structure} for the commutator of vector fields $[A,B]$. This is due to the fact that $[A,B]\in [\mathcal{R},\mathcal{R}]$, for $A,B \in \mathcal{R}$ and $ [\mathcal{R},\mathcal{R}]\neq \mathcal{R}$, see Remark \ref{r:spaceR}; further, the flat operator $\flat$ is only defined for matrices in $\mathcal{R}$.
To obtain a discrete vector $W\in\mathcal{R}$ approximating the commutator $\Big[\frac{\delta C}{\delta \mathbf{m}},\mathbf{u}\Big]$
at the edge midpoint, we will use the standard operators from Eq. \eqref{eq:discOps}. Then, we can compute $P(\mathcal{L}_A \left(hW^\flat\right))_{ij}$ using \cite[Lemma 3.1]{baue17a}.

\paragraph{Discrete commutator.} \label{sec:DiscreteComm}  Let $U_{ij}=\mathbf{u}_{e_{ij}}\cdot \mathbf{n}_{e_{ij}}$ be the edge normal for a vector field $\mathbf{u}$ at edge $e_{ij}$ and $V_{ij}$ for a vector field $\mathbf{v}$ respectively. 
The Lie bracket for vector fields $\mathbf{u}$ and $\mathbf{v}$ is given below and can be rewritten using a standard vector calculus identity, giving
$$
[\mathbf{u},\mathbf{v}]= \mathbf{u} \cdot \nabla \mathbf{v}- \mathbf{v}\cdot \nabla \mathbf{u}
=\mathbf{u}  \,\operatorname{div} \mathbf{v}   -\mathbf{v}  \,\operatorname{div} \mathbf{u}  -\nabla \times (\mathbf{u}\times \mathbf{v}).
$$
We discretize $\mathbf{u}  \,\operatorname{div} \mathbf{v} $ and $\mathbf{v} \,\operatorname{div} \mathbf{u}$ using the discrete divergence on each triangle (Eq.~\eqref{eq:discOps}), averaging over the two adjacent triangles to obtain an edge value,
%
\begin{align*}
	\Big(\mathbf{u}  \,\operatorname{div} \mathbf{v} \Big)_{ij}&=U_{ij}\Big(\frac{\Div(V)_i+\Div(V)_j}{2}\Big),
	\\
	\Big(\mathbf{v}  \,\operatorname{div} \mathbf{u} \Big)_{ij}&=V_{ij}\Big(\frac{\Div(U)_i+\Div(U)_j}{2}\Big).
\end{align*}
%
Then, to obtain a discrete version of $\nabla \times (\mathbf{u}\times \mathbf{v})$,  we use the following procedure:
\begin{itemize}
	\item Reconstruct the full vector fields $\mathbf{u}_\zeta$ and $\mathbf{v}_\zeta$ at the dual cell centres  
	from the normal values $U_{ij}$ and $V_{ij}$. We use the reconstruction in the interior of each triangle proposed by Ref. \cite{perot2006mimetic} and map it to the dual cell:
	\begin{align*}	 		 
		\mathbf{u}_\zeta &= \sum_{i\in N(\zeta)} \frac{|\zeta\cap T_i|}{|\zeta|} \mathbf{u}_i,
		\\
		\mathbf{u}_i &= \frac{1}{\Omega_{ii}} \sum_{k\in\{j,i_-,i_+\}}|e_{ik}|(\mathbf{x}_{e_{ik}}-\mathbf{x}_{T_i})U_{ij},	 
		\\
		\mathbf{v}_\zeta &= \sum_{i\in N(\zeta)} \frac{|\zeta\cap T_i|}{|\zeta|} \mathbf{v}_i,
		\\
		\mathbf{v}_i &= \frac{1}{\Omega_{ii}} \sum_{k\in\{j,i_-,i_+\}}|e_{ik}|(\mathbf{x}_{e_{ik}}-\mathbf{x}_{T_i})V_{ij}.	 
	\end{align*}	 
	%
	\item Compute the cross product  $c_\zeta=(\mathbf{u}_\zeta \times \mathbf{v}_\zeta)\cdot \mathbf{k}_\zeta$, where $\mathbf{k}$ is the unit vector that points in the local vertical direction.  
	\item Obtain $\Big(\nabla \times (\mathbf{u}\times \mathbf{v})\Big)_{ij}=\text{Grad}_t ~c_{\zeta}$. Since $c_\zeta$ is located at the dual cell centres and the resulting value after taking the curl should be an edge normal value (tangential for the dual grid cells), we use $\text{Grad}_t$, which is the adjoint curl
	(see Remark \ref{r:curlmap}).
\end{itemize}

We obtain $W$ approximating $\Big[\frac{\delta C}{\delta \mathbf{m}},\mathbf{u}\Big]$, as $W_{ij}=\frac{|e_{ij}|}{2\Omega_{ii}}\widetilde{W}_{ij}$ and $W_{ii}=(\text{Div}~\widetilde{W})_i$, where
%
\begin{align*}
	\widetilde{W}_{ij}&=\left(\frac{\delta C}{\delta M}\right)_{ij}\Big(\frac{\Div(V)_i+\Div(V)_j}{2}\Big)
	\\
	&~~~~ -V_{ij}\Big(\frac{\Div\left(\frac{\delta C}{\delta M}\right)_i+\Div\left(\frac{\delta C}{\delta M}\right)_j}{2}\Big)
	\\
	&~~~~ -\text{Grad}_t~\Big((\frac{\delta C}{\delta \mathbf{m}}_{\zeta}\times \mathbf{u}_\zeta)\cdot \mathbf{k}_\zeta\Big)_{ij}.
\end{align*}
%
In Ref. \cite[Lemma 3.1]{baue17a}, the discrete projection of the Lie derivative $ \mathcal{L}_A$ is given by
\begin{align}\label{eq:discLieNormal}
	\begin{split}
		P\Big(\textcolor{black}{ \mathcal{L}_A}(h W^\flat)\Big)_{ij}&=
		(\Curl_{\zeta_{-}}\widetilde{W})
		\Big(   \frac{|\zeta_{-} \cap T_i|}{2 \Omega_{ii}  } \overline{h}_{ji_-}  |e_{ii_-}|   V_{ii_-} \\
		& \hspace{6em}  +  \frac{|\zeta_{-}  \cap T_j|}{2 \Omega_{jj}  } \overline{h}_{ij_-}  |e_{jj_-}|   V_{jj_-} \Big) 
		\\& +
		(\Curl_{\zeta_{+}}\widetilde{W}) 
		\Big(  \frac{|\zeta_{+}  \cap T_i|}{2 \Omega_{ii} }\overline{h}_{ji_+} |e_{ii_+}|   V_{ii_+} \\
		& \hspace{6em}  + \frac{|\zeta_{+}  \cap T_j|}{2 \Omega_{jj}   } \overline{h}_{ij_+}   |e_{jj_+}|   V_{jj_+} \Big)
		\\&+
		\overline{h}_{ij} \Big( \sum_{k\in N(i)} \frac{|e_{ik}| |\tilde{e}_{ik}|}{\Omega_{ii}} V_{ik}\widetilde{W}_{ik} \\
		& \hspace{2em} -\sum_{k\in N(j)} \frac{|e_{jk}| |\tilde{e}_{jk}|}{\Omega_{jj}} V_{jk}\widetilde{W}_{jk}		
		\Big)
		\\&+
		\frac{\Div(V\overline{h})_i+\Div(V\overline{h})_j}{2}(2|\tilde{e}_{ij}|\widetilde{W}_{ij}),
	\end{split}
\end{align}
which we can evaluate using the discrete operators above once we know $\widetilde{W}$. \textcolor{black}{Moreover, since $W\in\mathcal{R}$, we use the definition of the flat operator \eqref{eq:discrFlat} to compute $W^\flat$ from $W$.}

\paragraph{Discrete enstrophy variational derivative.}
Analogously to Section \ref{sec:contEnstCas}, we compute the variational derivative of the approximation of the enstrophy Casimir and substitute it into Eq. \eqref{eq:momentumCas}. The discrete enstrophy Casimir is
\begin{align*}
	\mathcal{C}(M,h)&=\frac{1}{2} \sum_{\zeta } h_{\zeta} \Big(q(M,h)_\zeta\Big)^2 |\zeta|,	
	\\
	q(M,h)_\zeta &=\frac{(\text{Curl}~V) + f}{h_\zeta},
	\\
	h_\zeta &= \sum_{T_i \cap \zeta \neq \emptyset} \frac{|T_i\cap \zeta|}{|\zeta|} h_i,
\end{align*}
where $M=\frac{\delta \ell}{\delta A}$ and $f$ is the Coriolis parameter.

Then, computing the variational derivative (see Appendix \ref{A:vardercas} for details) we obtain
%
\begin{align*}
	\left( \frac{\delta \C}{\delta M}\right)_{ij} 
	&=
	\frac{q_{\zeta_+}-q_{\zeta_-}}{\Omega_{ii} h_{ij}}
	=
	-\frac{|e_{ij}|}{2\Omega_{ii}} \Big(2\frac{q_{\zeta_-}-q_{\zeta_+}}{|e_{ij}|} \frac{1}{h_{ij}}\Big)
	\\
	&=-\frac{|e_{ij}|}{2\Omega_{ii}} \frac{2~\text{Grad}_t ~q}{h_{ij}}.
\end{align*}
%
Substituting this into Eq. \eqref{eq:momentumCas} results in the discretized potential enstrophy dissipating SWE.
\begin{remark}
	The approximation of the enstrophy Casimir is not a Casimir of the discrete system. Therefore, we cannot directly prove that enstrophy is dissipated \textcolor{black}{(or does not grow)} for the semi-discrete scheme. However, the numerical results demonstrate that the numerical scheme indeed dissipates enstrophy.
\end{remark}

\subsection{Temporal discretization}
A temporal variational discretization can be obtained by following the discrete (in time) Euler--Poincar\'e--d'Alembert approach, see Refs. \cite{gawlik2011geometric,desbrun2014variational,GaGB2019}. This approach is based on the Cayley transform, a local approximation to the exponential map of the Lie group. In particular, the resulting scheme uses the Cayley transform in the update for the continuity equation and a Crank--Nicolson-type update for the momentum equation \textcolor{black}{given in \eqref{eq:projectionCas} with $ \theta =0$}. For the selective decay, the dissipation term is added to the Crank--Nicolson-type update. \textcolor{black}{Following Ref. \cite{baue17a}, we will use below the Crank--Nicolson-type time update directly on the momentum equation as reformulated in \eqref{eq:momentumCas}. This considerably simplifies the solution procedure without altering the behavior of the scheme.}

Based on the Cayley transformation, the continuity update equation  is then given by $h^{t+1}= \tau (\Delta t A^t)h^t$
for the time $t$ and a time step size $\Delta t$, where the action of $\tau$ can be represented by solving
\begin{equation}\label{equ_D} 
	\big(I - \frac{1}{2} \Delta t A^t\big) h^{t+1} = \big(I + \frac{1}{2} \Delta t A^t\big) h^{t},
\end{equation}
with $I$ being the identity matrix. 
Then, we use the following fixed-point iteration to approximately solve the discrete momentum equation:
\begin{enumerate}
	\item Start loop over $k$ with initial guess as solution at time $t$: $ V^{*}_{k=0} = V^{t}$;
	\item Calculate updated velocity $V_{k+1}^{*}$ from the explicit equation:
	\begin{align*}
		\frac{V^{*}_{k+1} -V^{t}}{\Delta t}  &=-\frac{\operatorname{Adv}_{\rm }(V^*_k,h^{t+1})+\operatorname{Adv}_{\rm }(V^{t  },h^{t  })}{2} 
		\\&~~~-\frac{\operatorname{K}_{\rm }(V^*_k )+\operatorname{K}_{\rm }(V^{t  } )}{2}-\operatorname{G}(h^{t+1})
		\\
		&~~~-\theta\frac{\operatorname{L}_{\rm }(V^*_k,h^{t+1},\frac{\delta C}{\delta M}^t)+\operatorname{L}_{\rm }(V^{t  },h^{t  },\frac{\delta C}{\delta M}^t)}{2};		
	\end{align*}
	\item Stop loop over $k$ if $||V^{*}_{k+1} - V^{*}_k|| < \epsilon$ for a small positive $\epsilon$, take $V^{t+1} = V^{*}_{k+1}$.
\end{enumerate}

For more details, we refer the reader to Refs. \cite{baue17a,brecht2019variational}.

\begin{remark}\label{r:temp}\textcolor{black}{For this scheme}, it has been observed that the \textcolor{black}{fully discrete} temporal integrator does not conserve energy at the level of machine precision but, rather, the energy error fluctuates around a long term mean.  Thus, while the energy is conserved by the semi-discrete Casimir dissipative equations, independent of the discretization of the commutator (see Proposition \ref{prop:compEnergy}), this does not guarantee that the energy will be conserved after temporal discretization.  We observe this in the numerical results that follow, but ascribe the small energy growth seen there to errors from this temporal discretization.  Attenuating these errors (or developing a tractable fully variational time integrator) is an open question for future research. \textcolor{black}{Note however that for the purpose of this work which focuses on the Casimir dissipation mechanism as a means to remove the small scale noise without interfering with the energy behavior, the given order of energy conservation of the temporal integrator is sufficient. Also, the use of an integrator that conserves energy at machine precision, does not preclude the need to incorporate the Casimir dissipation term to better preserves the coherent structures of the solution.}
\end{remark}

\subsection{Biharmonic dissipation}

\textcolor{black}{For the numerical simulations, we will compare the selective Casimir dissipation proposed here against simulations using no dissipation and those using a standard dissipation method. A common approach \textcolor{black}{that is used} to remove small scale noise and improve the stability of the scheme  is to apply a linear fourth-order diffusion (biharmonic dissipation)
	to the  velocity  field, see, e.g., Refs. \cite{Ripodas2009,flye12a,ring02a}}. For instance, on $ \mathbb{R} ^2 $, \textcolor{black}{this yields}:
\begin{equation}\label{eq:RSWdiss}
	\partial_t \mathbf{u} + \left(\nabla \times (\mathbf{u}+\mathbf{r})\right)\times \mathbf{u}+\nabla \left(\frac 1 2 |\mathbf{u}|^2+g(h+\eta_b)\right)=-\nu \Delta^2\mathbf{u},
\end{equation}
where $\nu$ is the diffusion coefficient. 
We do not add any dissipation to the continuity equation, because it does not contain a turbulent mixing term. Also, adding dissipation \textcolor{black}{to the continuity equation} can break conservation of mass, see Ref. \cite{Ripodas2009}.

We  discretize the dissipation term in Eq. \eqref{eq:RSWdiss}, using the vector calculus identity for the vector Laplacian 
$$
\Delta \mathbf{u}=\nabla  \,\operatorname{div} \mathbf{u} -\nabla\times(\nabla\times \mathbf{u}).
$$
Then, using the discrete operators \eqref{eq:discOps}, we obtain
\begin{align*}
	{\rm lap}(V)_{ij}&=\text{Grad}_n(\text{Div}~V)_{ij}-\text{Grad}_t(\text{Curl} ~V)_{ij}\\
	\nu \text{L}(V)_{ij}&=\nu ~{\rm lap}({\rm lap}(V))_{ij}. 	
\end{align*}
We obtain the discrete version of Eq. \eqref{eq:RSWdiss}
\begin{equation}\label{eq:momentumBi}
	\partial_t V_{ij}=-\text{Adv}(V,h)_{ij}-\text{K}(V)_{ij}-\text{G}(h)_{ij} -\nu \text{L}(V)_{ij}.
\end{equation}
The temporal discretization is the same as above, but with  $-\nu\text{L}(V)_{ij}$ instead of $\theta\operatorname{L}_{\rm }(V^{t  },h^{t  },\frac{\delta C}{\delta M}^t)$.
\textcolor{black}{Equation \eqref{eq:momentumBi} will be used below in our numerical results as a comparison with the selective Casimir dissipation approach that we propose in this paper.}

\section{Numerical results}\label{sec:numresults}
%
The numerical simulations in the plane are performed on a doubly periodic  rectangular domain $M  = [0, L_x] \times [0, L_y]$ with $L_x = 5000~\text{km}$ and $L_y = 4330~ \text{km}$. We consider an $f$-plane approximation with constant  Coriolis parameter $f$  set to $6.147\times 10^{-5}~\text{s}^{-1}$ and $g=9.81~\text{m}/\text{s}$. Unless otherwise noted, the simulations are performed using a resolution of $N=32768$ triangles. 
\textcolor{black}{
	For these simulations, we compute reference simulations on meshes with $524288$ triangles, using the scheme developed in Ref. \cite{mcra13Ay}, which conserves energy and uses the anticipated vorticity method for potential enstrophy dissipation.  
}

For the simulations on the sphere, we use an icosahedral grid, \textcolor{black}{as is widely adopted, see, for example, Refs. \cite{Ripodas2009,Aechtner15,heikes1995numericalP1}. We note, however, that this grid is optimized for the properties of its hexagonal and pentagonal dual cells, so that the numerical operators acting on these cells have a good convergence behaviour, at the expense of good convergence of similar operators on the triangular cells.} We set the Earth's radius $R=6.37122\times 10^6~\text{m}$, the Coriolis parameter to be $f=2\Omega\sin(\Theta)$, where $\Omega=7.292\times 10^{-5}~\text{s}^{-1}$, and $g=9.81~\text{m}/\text{s}$. Here, $\Theta$ is the latitude and $\Lambda$ the longitude. The simulations are performed using a resolution of $N=81920$ triangles and \textcolor{black}{a reference simulation using biharmonic dissipation is performed on a resolution of $N=327680$ triangles}.
\\~


We define the discrete total energy $H$, namely the Hamiltonian, and  the discrete  potential enstrophy $ \mathcal{C}  $:
\begin{align}
	H&=\sum_{T_i} \frac{g}{2}(h_i+(\eta_b)_i)^2\Omega_{ii}+\frac{1}{2}\Omega_{ii}h_i \sum_{k=j,i_-,i_+}\frac{|e_{ik}|~|\tilde{e}_{ik}|V_{ik}^2}{2}
	\\
	\mathcal{C} &=\frac{1}{2}\sum_{ \zeta }\frac{\text{Curl}(V(t))+f}{h_\zeta(t)}|\zeta|.
\end{align}

For each \textcolor{black}{test case}, we first choose the dissipation coefficient $\nu$ \textcolor{black}{for the biharmonic dissipation simulation}, and then empirically choose $\theta$ for the Casimir dissipation \textcolor{black}{simulation} so that the dissipation of potential enstrophy is qualitatively similar \textcolor{black}{between these two models}.
\subsection{Numerical analysis of the discrete commutator}

We first present a convergence study for the discrete commutator on both the plane and sphere. We define
\begin{align*}
	&\mathbf{u}=\begin{pmatrix}
		\sin(\frac{2\pi x}{L_x})\\0
	\end{pmatrix} \text{ and }\mathbf{v}=\begin{pmatrix}
		\cos(\frac{2\pi x}{L_x})\\0
	\end{pmatrix}, 
	\\
	&\text{ such that } [\mathbf{u},\mathbf{v}]=\begin{pmatrix}
		\frac{2\pi}{L_x}\\0
	\end{pmatrix}
\end{align*}
for the test case in the plane, and 
$$
\mathbf{u}=\begin{pmatrix}
y\\-x\\0
\end{pmatrix} \text{ and }\mathbf{v}=\begin{pmatrix}
0\\-z\\y
\end{pmatrix} \text{ such that } [\mathbf{u},\mathbf{v}]=\begin{pmatrix}
z\\0\\-x
\end{pmatrix}
$$
for the test case on the sphere. 
We approximate these vector fields with piecewise constant functions and follow  the algorithm in Section \ref{sec:DiscreteComm}
to discretize $[\mathbf{u},\mathbf{v}]$. Then, we compute the error between the approximation of the discrete commutator to the analytic field projected on the edge normal direction.  To estimate the numerical errors, we use the following definitions for the relative $L_2$ and $L_\infty$ error on edge values:
\begin{align*}
	L_2 &=\frac{\sqrt{\sum_{ij}|\mathfrak{e}_{ij}|\Big(u_n(e_{ij})-u_r(e_{ij})\Big)^2}}{\sqrt{\sum_{ij}|\mathfrak{e}_{ij}|u_r(e_{ij})^2}},
	\\
	L_\infty &= \frac{\max_{ij}|u_n(e_{ij})-u_r(e_{ij})|}{\max_{ij}|u_r(e_{ij})|},
\end{align*}
where $u_n(e_{ij})$ is the numerical solution defined at edge $e_{ij}$ and $u_r(e_{ij})$ is the analytical solution evaluated at the edge midpoint $e_{ij}$. Moreover $|\mathfrak{e}_{ij}|=\frac{1}{2}|e_{ij}|~|\tilde{e}_{ij}|$ is the area associated to an edge.

\begin{figure}
	\begin{center}
		\includegraphics[width=0.5\textwidth]{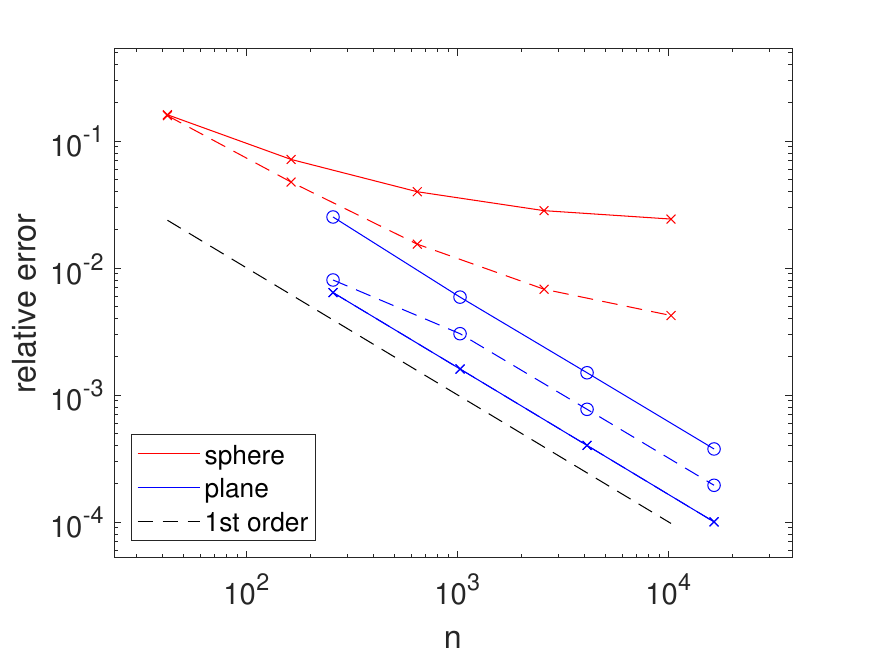}
	\end{center}
	\caption{Convergence of the commutator on the plane and sphere. Solid lines represent the $L_\infty$ error and dashed lines the $L_2$ error. In the plane, we present results using both a regular grid, with data denoted by x, and an irregular grid (with a central refinement region), with data denoted by o.}
	\label{fig:convCommutator}
\end{figure}

We observe that the resulting approximations on both a regular and irregular grid in the plane are first-order accurate, see Fig. \ref{fig:convCommutator}. On the sphere, the approximation is  less than first-order accurate.
\textcolor{black}{This is expected because the icosahedral grid is optimized for properties of the hexagonal and pentagonal dual cells and not those of the triangular cells used here. Moreover, Ref. \cite{heikes1995numerical} notes that, without the optimization of the grid for the triangular cells, the numerical operators converge at less than first order. Here in particular, we evaluate the discrete divergence and reconstruction of the vector fields on the non-optimized triangles. This low-order convergence of the divergence was also observed in Ref. \cite{brecht2019variational}. }

\subsection{Vortex interaction}
This test case consists of two counter-rotating vortices in the plane and is dominated by nonlinear processes. The two vortices are initially placed too far apart to merge. Thus, a key point in this simulation is that adding the Casimir dissipation does not change the evolution of the vortices.
\paragraph{Initial conditions.}
The initial height function for this example \cite{baue17a} is given by 
\begin{subequations}
	\begin{align}
		\begin{split}
			h \big( x, y, t=0 \big) &= H_0
			- H' \bigg( \exp \Big( -\frac{{x'_1}^2+{y'_1}^2}{2} \Big) 
			\\&~+ \exp \Big(-\frac{{x'_2}^2+{y'_2}^2}{2} \Big) -\frac{4\pi s_x s_y}{L_x L_y} \bigg),
		\end{split}
	\end{align}
	where $H_0=750\,m$, $H'=75\, m$, and the periodic extensions are given by
	\begin{align}
		\begin{split}
			x'_i &= \frac{L_x}{\pi s_x} \sin \big( \frac{\pi}{L_x} (x - x_{c_i}) \big)
			\\
			y'_i &= \frac{L_y}{\pi s_y} \sin \big( \frac{\pi}{L_y} (y - y_{c_i}) \big),\ \quad
			i = 1,2
		\end{split}
	\end{align}
	with the centres located at $(x_{c_1}, y_{c_1}) = 2/5\, (L_x, L_y)$, $(x_{c_2}, y_{c_2}) = 3/5\, (L_x, L_y)$ and $(s_x , s_y) = 3/40\, (L_x, L_y)$.
\end{subequations}
The discrete initial water depth on each triangle, $h_i$, is obtained by sampling  the analytical water depth at the cell center. Then, the initial condition for the velocity is given by the discrete geostrophic velocity,
$$
V_{ij}=-\frac{g}{f}\text{Grad}_t(h)_{ij}.
$$
In these simulations, we use \textcolor{black}{$\Delta t=0.00069\,$day} and dissipation parameters $\nu=1.2724\times 10^5~\text{km}^4/\text{day}$ and $\theta=2~\text{km}^4\text{day}^2$. \textcolor{black}{For the reference simulation, we use $\Delta t=0.005\,$day and set the APVM parameter to $\Delta t/200$.}

\begin{figure}
	\centering 
	\includegraphics[width=0.4\textwidth]{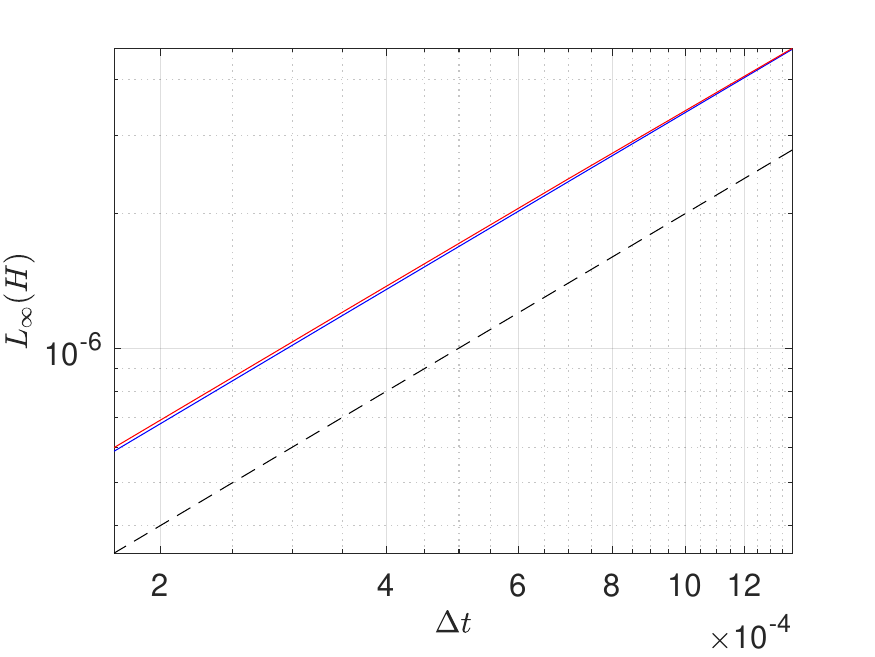}
	\caption{Convergence of the energy with respect to the time step size on a regular (blue) and irregular (red) grid. The dashed black line indicates first order.}\label{fig:enconv}
\end{figure}

\begin{figure*}
	\begin{center}
		\includegraphics[clip, trim=3cm 0cm 3cm 0cm,width=\textwidth]{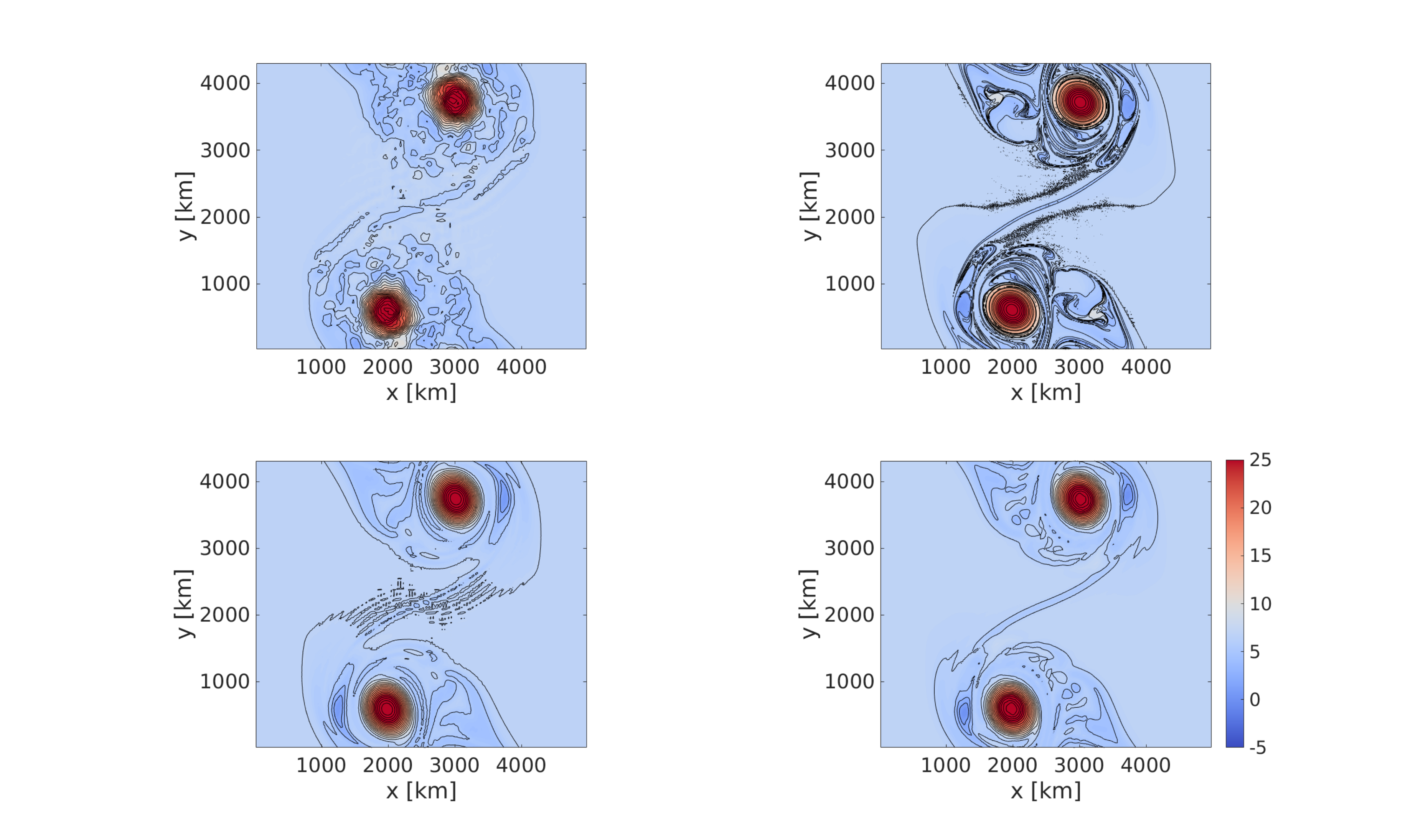}
	\end{center}
	\caption{Interacting vortices test case: Comparison of the potential vorticity for a simulation without dissipation (top left), the reference simulation (top right), a simulation with Casimir dissipation (bottom left) and a simulation with biharmonic dissipation (bottom right) after 10 days.  }\label{fig:compvtx}
\end{figure*}

We first integrate the initial conditions for two days for different
values of the time step, to analyze the convergence of the energy. In Fig. \ref{fig:enconv}, we observe that the energy converges with first-order accuracy.
Then, to analyze the effects of the Casimir dissipation, we integrate the initial conditions for 10 days and compare the relative potential vorticity field against a simulation with no dissipation, \textcolor{black}{the reference simulation,} and one with biharmonic dissipation, see Fig.  \ref{fig:compvtx}. All simulations behave similarly, with the cores of the two vortices being mutually repelled, due to nonlinear effects. We note that the simulation with no dissipation becomes noisy, while the two simulations with dissipation retain their accuracy. 
\textcolor{black}{When comparing to the reference simulation, the Casimir and biharmonic dissipation simulation behave similarly. However, the biharmonic dissipation simulation dissipates more small-scale motion.}

%

\begin{figure*}[h!]
	\begin{center}
		\includegraphics[width=\textwidth]{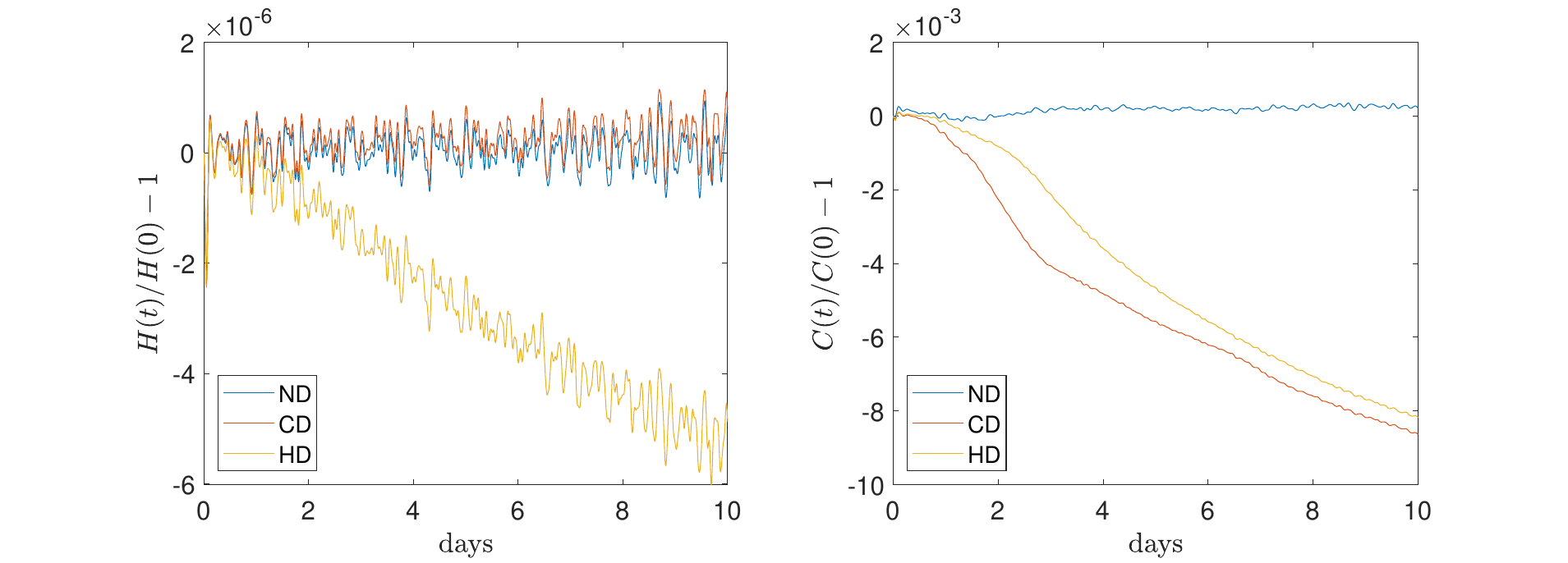}
	\end{center}
	\caption{Interacting vortices test case: Comparison of the relative errors in the energy (left) and potential enstrophy (right)  for a simulation without dissipation (blue), with Casimir dissipation (red) and with a biharmonic dissipation (yellow).}\label{fig:vtxDiagnostics}
\end{figure*}

The quantities of interest, total energy and potential enstrophy, are shown in Fig. \ref{fig:vtxDiagnostics}. We observe that the enstrophy is dissipated at the same rate for the simulations with biharmonic and Casimir dissipation, as expected with this choice of dissipation parameters. While the energy is dissipated in the simulation with biharmonic dissipation, conservation of energy for the simulation using Casimir dissipation is similar to that of the simulation with no dissipation. 
As noted in Remark \ref{r:temp} above, the temporal discretization used here is not completely energy conserving, leading to the oscillations seen at the left of Fig. \ref{fig:vtxDiagnostics}.

\subsection{Shear flow}
We next consider a shear flow test case in the quasi-geostrophic regime \cite{baue17a}, with strongly dominant nonlinear effects. The shear flow is initialized to an unstable equilibrium state so that, after a few days, the instability develops. This test case demonstrates that adding the Casimir dissipation does not change the development and growth of this instability.

\paragraph{Initial conditions.} The initial height for this example is given by 
\begin{align*}
	h(x,y,t=0)&=H_0-H'\frac{y''}{\sigma_y}e^{-\frac{{y'}^2}{2\sigma_y^2}+\frac 1 2}\Big( 1-\kappa\sin\Big(\frac{2\pi x'}{\lambda_x}\Big) \Big),
\end{align*}
where
%
\begin{align*}
	x'&= \frac{x}{L_x},
	\\
	y'&=\frac 1 \pi \sin\Big(\frac{\pi}{Ly}\Big(y-\frac{L_y}{2}\Big)\Big),\qquad 
	\\
	y''&=\frac{1}{2\pi}\sin\Big(\frac{2\pi}{L_y}\Big(y-\frac{L_y}{2}\Big)\Big),
\end{align*}
%
with parameters $\lambda_x=\frac 1 2, \sigma_y=\frac 1 {12}, \kappa=0.1, H_0=1.076~\text{km}$ and $H'=0.03~\text{km}$. Again,  the velocity field is initialized to be the discrete geostrophic velocity,
$$
V_{ij}=-\frac{g}{f}\text{Grad}_t(h)_{ij}.
$$
Here, \textcolor{black}{ $\Delta t=0.010\,$day and} the dissipation parameters are chosen as $\nu=3.7145\times 10^5~\text{km}^4/\text{day}$ and $\theta=2~\text{km}^4\text{day}^2$. \textcolor{black}{For the reference simulation, we use $\Delta t=0.025\,$day and set the APVM parameter to $\frac{\Delta t}{20}$.}

\begin{figure*}
	\begin{center}
		\includegraphics[clip,trim=3cm 0cm 3cm 0cm,width=\textwidth]{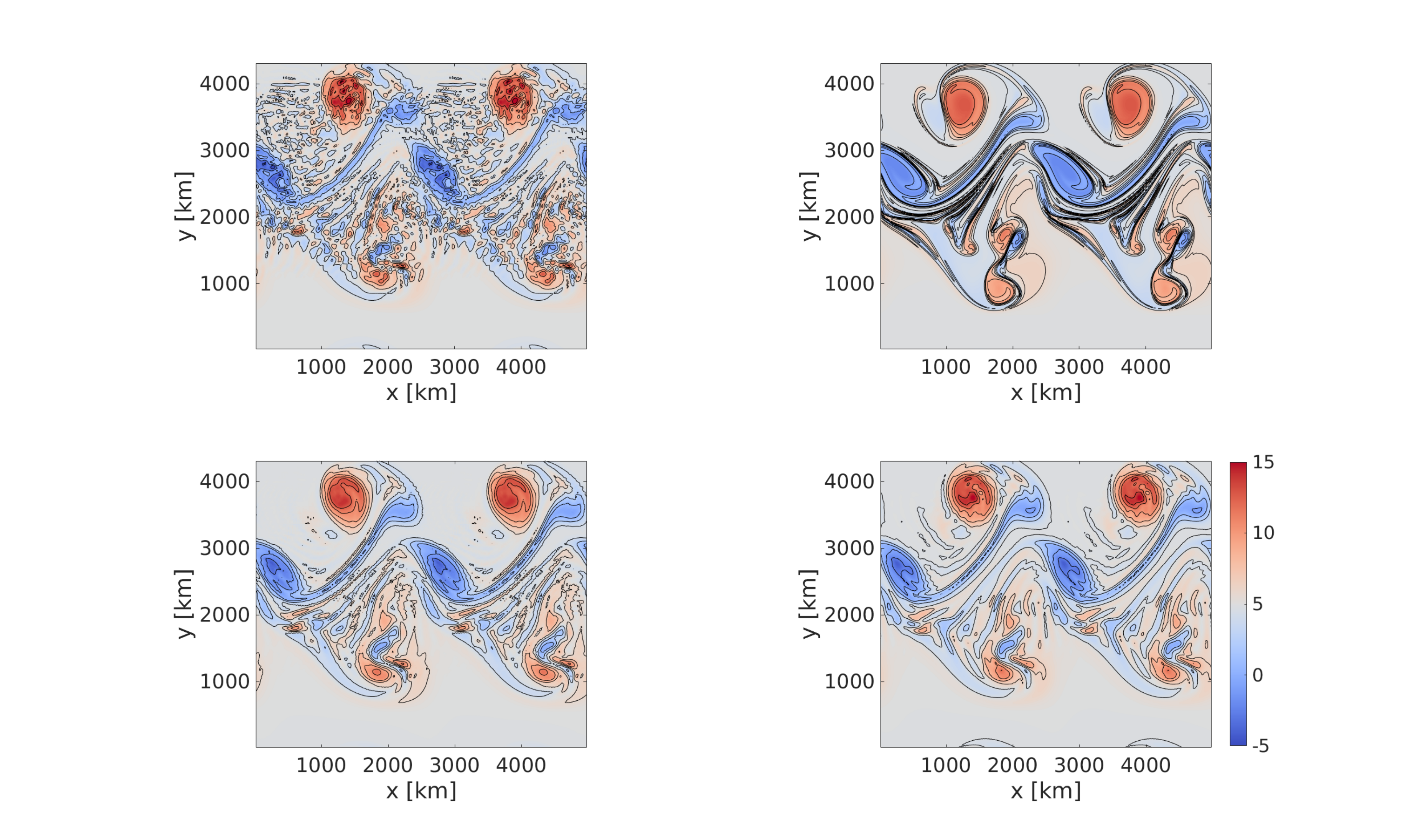}
	\end{center}
	\caption{Shear flow test case: Comparison of the  potential vorticity for a simulation without dissipation (top left), the reference simulation (top right), a simulation with Casimir dissipation (bottom left) and a simulation with biharmonic dissipation (bottom right) after 10 days.}\label{fig:shearComp}
\end{figure*}

We integrate the initial conditions for 10 days. The instability develops in the first three days, then the flow evolves into pairs of counter-rotating vortices. The filaments between the vortices become thinner until they can no longer be resolved by the spatial resolution of the mesh. This causes a noisy pattern in the vorticity field at day ten for the simulation without any dissipation, see Fig. \ref{fig:shearComp}. In contrast, the simulations with Casimir and biharmonic dissipation are much less polluted. \textcolor{black}{Moreover, the Casimir dissipation appears to preserve more of the fine-scale structure seen in the reference simulation than does the simulation with biharmonic dissipation.}

\begin{figure*}
	\begin{center}
		\includegraphics[width=\textwidth]{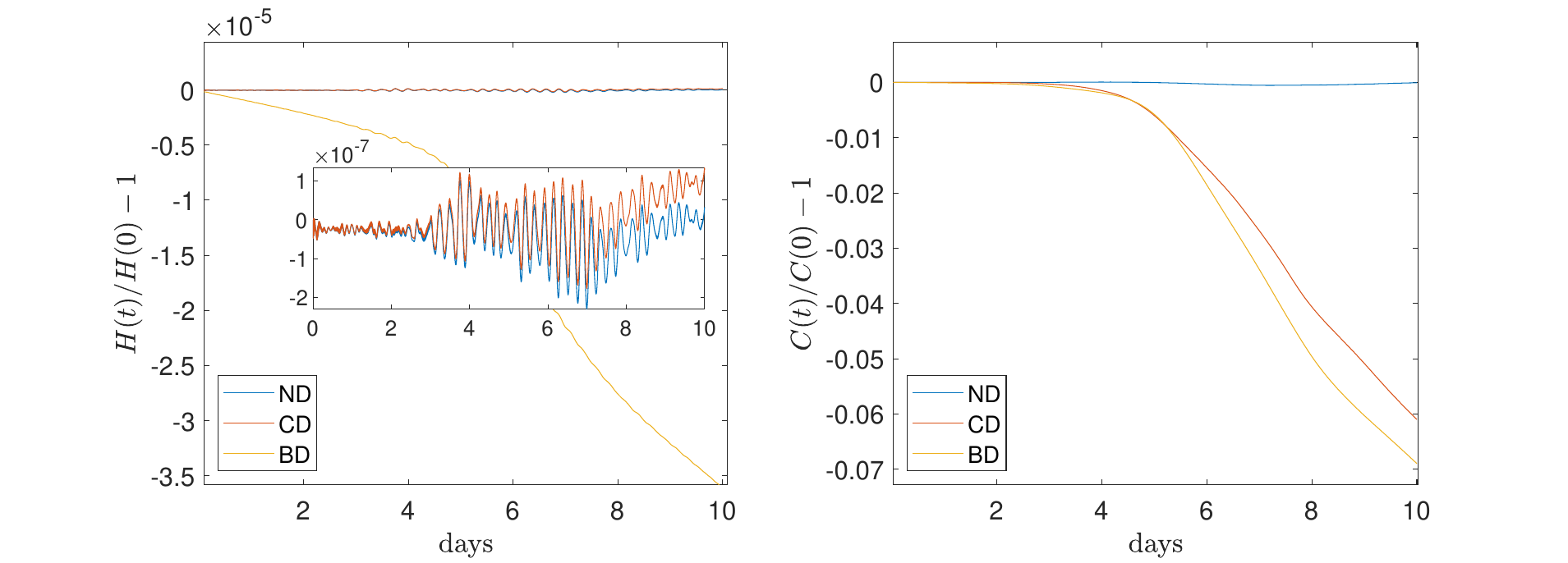}	
	\end{center}
	\caption{Shear flow test case: Comparison of the relative errors in the energy (left) and potential enstrophy (right) for a simulation without dissipation (blue), with Casimir dissipation (red) and with biharmonic dissipation (yellow).}
	\label{fig:shrDiagnostics}
\end{figure*}

The quantities of interest for this simulation are shown in Fig. \ref{fig:shrDiagnostics}.  Again, we observe the similar dissipation rate of the potential enstrophy for the Casimir and biharmonic dissipation, by construction. The simulation with no dissipation and the Casimir dissipative simulation have a similar conservation of energy.  In contrast, the simulation with biharmonic dissipation has a loss of energy about 100 times greater.

\begin{figure*}
	\begin{center}
		\includegraphics[width=\textwidth]{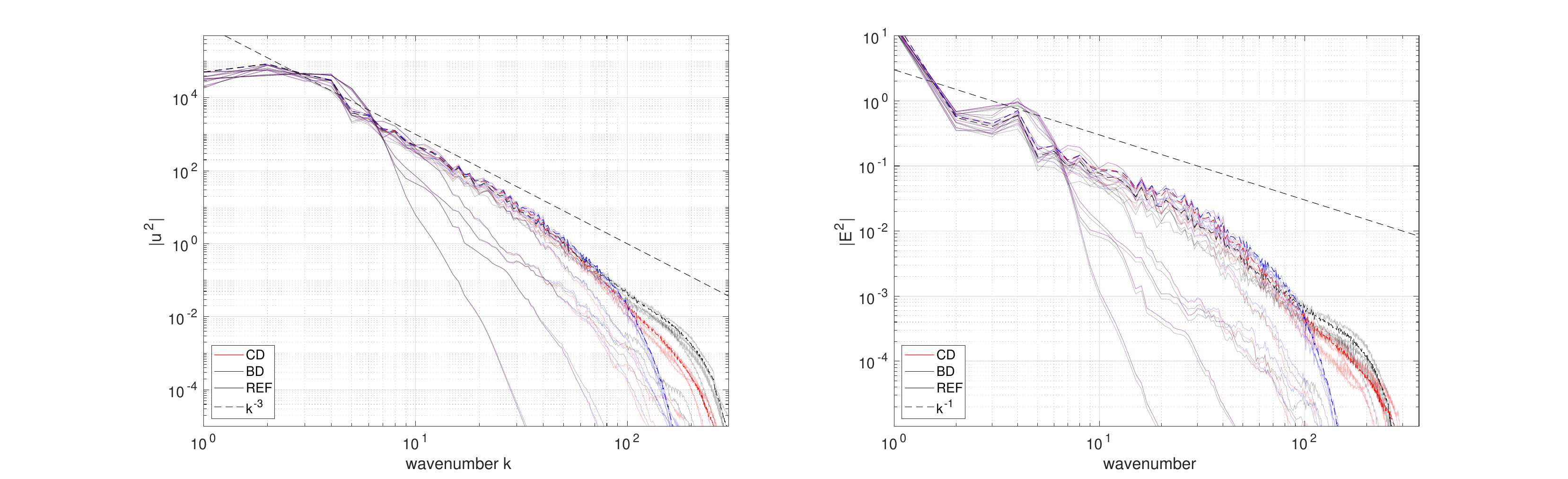}		
	\end{center}
	\caption{Shear flow test case: Comparison of the kinetic energy (left) and potential enstrophy (right) spectra for the reference simulation (black) and simulations with Casimir dissipation (red) and  biharmonic dissipation (blue). The spectra are obtained from simulations with a resolution of $N=524288$ triangles. The spectra are shown for days 1 to 10 of the simulation. The dashed red and blue lines show the averages of the spectra from days 6 to 10.}
	\label{fig:shrSpc}
\end{figure*}

Fig. \ref{fig:shrSpc} shows the kinetic energy and potential enstrophy spectra for simulations on refined spatial meshes with 524288 triangles.  Expected scaling laws for these spectra are discussed in Refs. \cite{ring02a,chen11Ay}. Both dissipative simulations follow the expected $k^{-1}$ power law for the enstrophy and $k^{-3}$ power law for the kinetic energy over a significant region of the resolved wavenumbers. However, we note that using Casimir dissipation results in better resolution of the spectra over the small scales (higher wavenumbers) in comparison with the biharmonic dissipation.  As expected, the biharmonic dissipation results in much faster dissipation over small scales in both the energy and enstrophy.
\subsection{Flow over an isolated mountain}
As a final example, 
we consider the flow over a conically-shaped mountain on the sphere, as proposed in Ref. \cite{wil92}. The initially balanced flow runs over the mountain, which initiates turbulence. The flow stays turbulent for a long period of time. 

\paragraph{Initial conditions.}
The discrete initial velocity and height fields are given in spherical coordinates as
\begin{align*}
	V_{ij}&=u_0(\cos(\Theta),0)^\top\cdot \mathbf{n}_{ij}
	&& u_0=20~\text{m}/\text{s}
	\\
	h_i&= h_0-\frac{1}{g}(R\Omega u_0+u_0^2/2)\cos(\Theta)
	&& h_0=5960~\text{m}
\end{align*} 
The conically shaped bottom topography is given by
%
\begin{align*}
	\eta_b(\Lambda,\Theta)&=2000(1-9r/\pi),
	\\
	r^2&=\min\Big( (\pi/9)^2, (\Lambda-\Lambda_c)^2+(\Theta-\Theta_c)^2\Big),
\end{align*}
%
where $\Lambda_c=3\pi/2$ and $\Theta_c=\pi/6$.

\textcolor{black}{Here, $\Delta t=100\,$s and the dissipation parameters are chosen as $\nu=1.9508\times 10^{14}~\text{km}^4/\text{day}$ and $\theta=-1\times 10^{20}~\text{km}^4\text{day}^2$. For the reference simulation, we use $\Delta t=25\,$s and $\nu=1.2191\times 10^{13}~\text{km}^4/\text{day}$.}

\begin{figure*}
	\begin{center}
		\includegraphics[clip, trim=3cm 2cm 3cm 2cm,width=\textwidth]{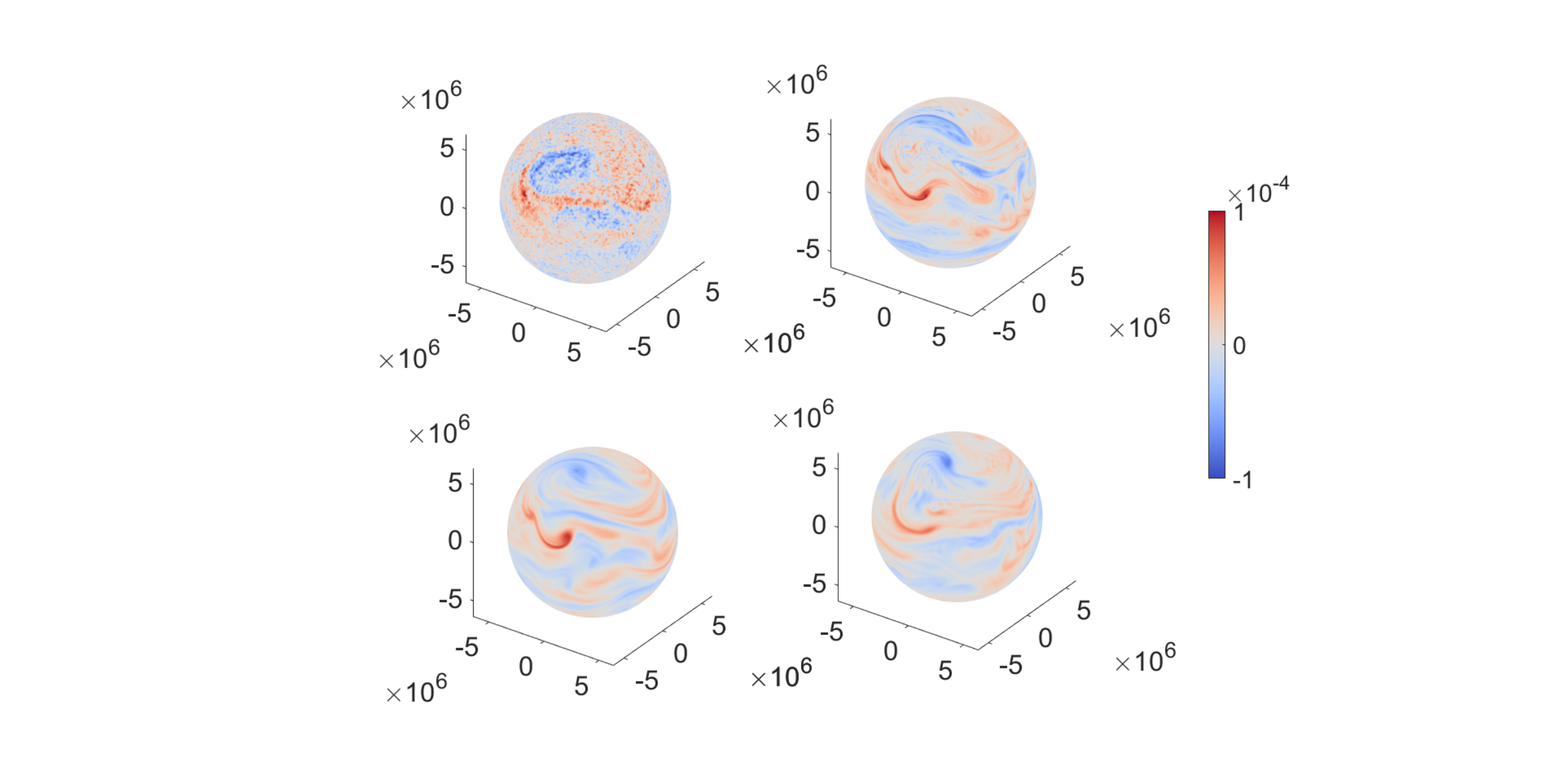}
	\end{center}	
	\caption{Flow over a mountain test case: Comparison of the relative vorticity for a simulation without dissipation (top left), the reference simulation (top right), a simulation with Casimir dissipation (bottom left), and a simulation with biharmonic dissipation (bottom right) after 100 days. }
	\label{fig:sphereFilds}
\end{figure*}

\textcolor{black}{We integrate the initial conditions for 100 days and compare a simulation without dissipation, a reference simulation, simulations with Casimir dissipation and biharmonic dissipation, see Fig. \ref{fig:sphereFilds}. The simulation without any stabilization becomes noisy, while the stabilized schemes  produce coherent structures in the vorticity field. When comparing the simulations with dissipation to the reference simulation, we observe that the fields are different. However, the Casimir dissipation simulation seems to match the vorticity field of the reference simulation better than the biharmonic dissipation simulation. This is because the energy is not dissipated away and more small scales are resolved.}

\begin{figure*}
	\begin{center}
		\includegraphics[width=\textwidth]{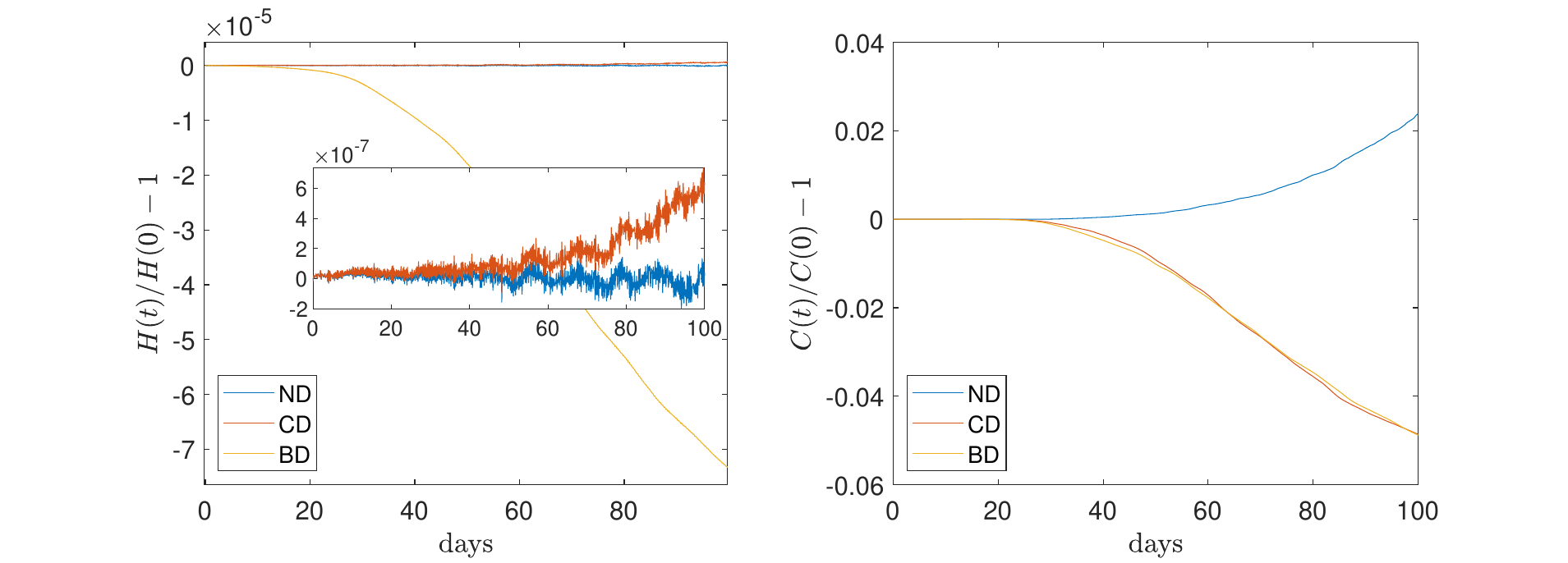}
	\end{center}
	\caption{Flow over mountain test case: Comparison of the relative errors in the energy (left) and potential enstrophy (right) for a simulation without dissipation (blue), with Casimir dissipation (red), and with biharmonic dissipation (yellow).}
	\label{fig:diagsphere}
\end{figure*}

In Fig. \ref{fig:diagsphere}, we show the quantities of interest. The simulation without dissipation shows an increase in potential enstrophy, which is related to the noisy vorticity field. The dissipative schemes, as expected, dissipate potential enstrophy at the same rate.
The error in the energy of the simulation with the Casimir dissipation stays on the same order as the simulation without any dissipation, while the simulation using biharmonic dissipation has an energy loss. \textcolor{black}{We note that there is a small, but consistent, gain in the energy for the simulation using Casimir dissipation. As noted above, the development of a fully conservative integrator in this setting is non-trivial, but a key question for future research.}


\section{Conclusions}\label{sec:conclusion}
The development of high-fidelity numerical simulation tools for weather and climate prediction is limited by the competing goals of achieving energy conservation while preserving long-term stability of the time-integration scheme (see e.g. Ref. \cite{sado85Ay}).  To address this challenge, we consider a potential-enstrophy dissipation scheme that conserves energy, building on existing variational integrators for the rotating shallow water equations \cite{baue17a,brecht2019variational}. In particular, the scheme combines the variational discretization framework of Ref. \cite{pavlov2011structure} with the selective decay proposed in  Ref. \cite{gay2013selective}.  The resulting semi-discrete scheme is shown to conserve energy, suggesting this is a viable framework for long-term climate simulations. 

Numerical results are presented comparing the new scheme with the variational integrator without dissipation and with a standard dissipation approach using a biharmonic eddy viscosity term. These simulations are carried out on both the $f$-plane and sphere, and we observe  that the simulations with no dissipation becomes noisy, in contrast to the simulations with dissipation. When analyzing the conservation properties, we find that the enstrophy dissipating scheme conserves the energy \textcolor{black}{to the same order as} the scheme without dissipation, while the biharmonic dissipation leads to a substantial loss of energy. Additionally, by computing energy and enstrophy spectra, we see that simulations using enstrophy dissipation better resolve small-scale features than those using biharmonic dissipation.
\textcolor{black}{In particular, the simulation on the sphere demonstrates the benefit of better-resolving these small scales, resulting in a simulation that more closely resembles that of the reference solution computed at higher resolution.}

\textcolor{black}{As shown above, the Casimir dissipation scheme stabilizes the variational integrator for a longer period of time, but neither scheme leads to a fully energy-conservative method after the temporal discretization. Thus, a key next step in this research is the development of a fully energy conserving variational temporal discretization for cubic Lagrangians, as considered here.}
Further natural extensions of this work are to more realistic models for geophysical flows.  In particular,
the primitive equations are a common first step in developing accurate simulations of climate dynamics.  
Since the framework in Ref. \cite{gay2013selective} also applies to 3D flows, combining the variational discretization framework with Casimir selective decay would lead to a discretization methodology for the primitive equations that would enable stabilized long-term integration schemes.
\textcolor{black}{Moreover, the method can be extended to include boundary conditions on non-periodic domains, as needed to apply the selective Casimir dissipation scheme in areas such as ocean modelling.}

\section{Acknowledgment}
This research was undertaken, in part, thanks to funding
from the Canada Research Chairs program, the NSERC
Discovery Grant program, and the InnovateNL CRC Leverage R\&D program.  WB would like to acknowledge funding from NERC NE/R008795/1. FGB was supported by the project ANR-14-CE23-0002-01.
\appendix
\section{Detailed computations}
\subsection{Continuous functional Casimir derivative}\label{sec:VarDerEnstCas}

The variational derivative of the Casimir on a doubly periodic planar domain is computed as follows
\begin{align*}
	\int_M \frac{\delta C}{\delta \mathbf{m} } \cdot \delta \mathbf{m} \,{\rm d} \mathbf{x}
	&= \left. \frac{d}{d\varepsilon}\right|_{\varepsilon=0}  C( \mathbf{m} +\varepsilon \delta \mathbf{m} , h)
	\\&=\left. \frac{d}{d\varepsilon}\right|_{\varepsilon=0} \frac{1}{2}\int_M h q( \mathbf{m} +\varepsilon \delta \mathbf{m} , h)^2 {\rm d} \mathbf{x}
	\\&=
	\textcolor{black}{
		\left. \frac{d}{d\varepsilon}\right|_{\varepsilon=0} \frac{1}{2}\int_M 
		\frac{\big(\mathbf{z}\cdot \text{curl}(\frac{ \mathbf{m} +\varepsilon \delta \mathbf{m} }{h})\big)^2}{h}{\rm d} \mathbf{x}}
	\\&=
	\textcolor{black}{
		\int_M \frac{\mathbf{z}\cdot \text{curl}(\frac{ \mathbf{m} }{h})}{h}\left(\mathbf{z}\cdot \text{curl}\left(\frac{\delta \mathbf{m} }{h}\right)\right) {\rm d} \mathbf{x}
	}
	\\&=\int_M q\mathbf{z} \cdot \text{curl}\left(\frac{\delta \mathbf{m} }{h}\right) {\rm d} \mathbf{x}
	\\&=\int_M \text{curl}(q\mathbf{z})\cdot \frac{\delta \mathbf{m} }{h} {\rm d} \mathbf{x}. 
\end{align*}
In the final step, we use the identity $\nabla \cdot (A\times B)=(\nabla\times A)\cdot B-A\cdot (\nabla \times B)$, noting that $\int_{M} \operatorname{div} \mathbf{u} \, {\rm d} \mathbf{x}  = \int_{\partial M}  \mathbf{u}  \cdot \mathbf{n}\, {\rm d} S=0$, since $M$ is doubly periodic and, thus, has no boundary.
This gives
\begin{equation}\label{eq:contDC}
	\frac{\delta C}{\delta \mathbf{m} } = \frac{1}{h} \operatorname{curl} (q \mathbf{z}) = - \mathbf{z} \times  ( \nabla q) / h
	=(\partial_y q, -\partial_x q,0)^\top/ h,
\end{equation}
where we first use the identity $\nabla \times (\psi A)=\psi (\nabla \times A) + \nabla \times A $ and, then, use the fact that $\mathbf{z}$ is the canonical unit vector in the $z$-direction.

Similarly, when $M$ is a two-dimensional Riemannian manifold, we compute
\begin{align*}
	\int_M \frac{\delta C}{\delta \mathbf{m} } \cdot \delta \mathbf{m} \,{\rm d}\sigma 
	&=\left. \frac{d}{d\varepsilon}\right|_{\varepsilon=0} \frac{1}{2}\int_M h q( \mathbf{m} +\varepsilon \delta \mathbf{m} , h)^2 {\rm d} \sigma 
	\\&= \int_M h q  \left. \frac{d}{d\varepsilon}\right|_{\varepsilon=0} q( \mathbf{m} +\varepsilon \delta \mathbf{m} , h){\rm d} \sigma \\
	\\&=\int_M  q \mathbf{d} \frac{ \delta \mathbf{m} }{h} {\rm d} \sigma 
	\\&= \int_M \mathbf{d} \left( \frac{ \delta \mathbf{m} }{h} q \right) + \int_M \frac{ \delta \mathbf{m} }{h} \wedge \mathbf{d} q
	\\&= - \int_M \frac{ \delta \mathbf{m} }{h} \cdot ( \star \mathbf{d} q)^\sharp {\rm d} \sigma,
\end{align*}
which gives $\frac{\delta C}{\delta \m} = -\frac{1}{h} (\star \mathbf{d}  q)^\sharp$. In the computation above, we have used Stokes' theorem on $M$, $\int_M \mathbf{d} \alpha =\int_{ \partial M} \alpha =0$ (since $ \partial M=\varnothing$), and the identity $\alpha \wedge \star \beta = (\alpha \cdot \beta ^\sharp) {\rm d} \sigma $, for one-forms $ \alpha , \beta $ on $M$, with $\star$ and $\sharp$ denoting the Hodge star and sharp operators associated with the Riemannian metric.

\subsection{Discrete functional Casimir derivative}\label{A:vardercas}

The semi-discrete variational derivative of the enstrophy is given by

	\begin{align*}
		\left\langle \frac{\delta \C}{\delta M}, \delta M \right\rangle_1  
		&= \left. \frac{d}{d\varepsilon}\right|_{\varepsilon=0} \C(M+ \varepsilon\delta M, h)
		\\
		\text{Tr}\Big(\frac{\delta \C}{\delta M}^\top\Omega \delta M\Big) 
		&=\left. \frac{d}{d\varepsilon}\right|_{\varepsilon=0} \frac{1}{2}\sum_{\zeta} \frac{1}{h{\zeta}} \left( \sum_{ \tilde{e}_{nm} \in \partial \zeta } \frac{M_{nm}+\varepsilon \delta M_{nm}}{| \zeta|h_{nm}} \right)^2| \zeta| 
		\\
		\sum_i\sum_j\frac{\delta \C}{\delta M}_{ij}^\top(\Omega \delta M)_{ji} 
		&= \sum_{\zeta} \frac{1}{h_\zeta} \left( \sum_{\tilde{e}_{nm} \in \partial \zeta} \frac{M_{nm}}{| \zeta|h_{nm}}
		\cdot \sum_{ h_{nm} \in \partial \zeta } \frac{\delta M_{nm}}{| \zeta |h_{nm}}
		\right)| \zeta| 
		\\
		\sum_{ij}\frac{\delta \C}{\delta M}_{ij}\Omega_{ii} \delta M_{ij} 
		&= \sum_{\zeta} \frac{1}{h_\zeta} \left( \sum_{ \tilde{e}_{nm} \in \partial \zeta} \frac{M_{nm}}{| \zeta|h_{nm}}
		\cdot \sum_{ h_{nm} \in \partial \zeta } \frac{\delta M_{nm}}{h_{nm}}
		\right)
		\\		
		\sum_{i}\sum_{j\in N(i)}\frac{\delta \C}{\delta M}_{ij}\Omega_{ii} \delta M_{ij} 
		&=\sum_{ \zeta } q_\zeta \sum_{ h_{nm} \in \partial \zeta } \frac{\delta M_{nm}}{h_{nm}}
		\\&=\sum_{ij} \frac{\delta M_{ij}}{h_{ij}} q_{\zeta_+} + \frac{\delta M_{ji}}{h_{ji}} q_{\zeta_-}
		\\&=
		\sum_{ij} \frac{\delta M_{ij}}{h_{ij}} (q_{\zeta_+} - q_{\zeta_-}).
	\end{align*}

In the second-to-last step, we use the property that each edge $e_{ij}$ has 2 neighboring vertices, denoted by $\zeta_+$ and $\zeta_-$.  In the last step, we use the fact that the matrix $M$ is anti-symmetric, as is $\delta M$, while $h$ is symmetric.

\bibliography{biblo,bihlo,bwerner1,bwerner}

\begin{thebibliography}{41}
\providecommand{\natexlab}[1]{#1}
\providecommand{\url}[1]{\texttt{#1}}
\expandafter\ifx\csname urlstyle\endcsname\relax
  \providecommand{\doi}[1]{doi: #1}\else
  \providecommand{\doi}{doi: \begingroup \urlstyle{rm}\Url}\fi

\bibitem[Aechtner et~al.(2015)Aechtner, Kevlahan, and Dubos]{Aechtner15}
M.~Aechtner, N.~K.-R. Kevlahan, and T.~Dubos.
\newblock A conservative adaptive wavelet method for the shallow-water
  equations on the sphere.
\newblock \emph{Quarterly Journal of the Royal Meteorological Society},
  141\penalty0 (690):\penalty0 1712--1726, 2015.

\bibitem[Arakawa and Hsu(1990)]{arakawa1990energy}
A.~Arakawa and Y.J.G. Hsu.
\newblock Energy conserving and potential-enstrophy dissipating schemes for the
  shallow water equations.
\newblock \emph{Monthly Weather Review}, 118\penalty0 (10):\penalty0
  1960--1969, 1990.

\bibitem[Bauer and Gay-Balmaz(2019)]{baue17a}
W.~Bauer and F.~Gay-Balmaz.
\newblock Towards a geometric variational discretization of compressible
  fluids: the rotating shallow water equations.
\newblock \emph{Journal of Computational Dynamics}, 6\penalty0 (1):\penalty0
  1--37, 2019.

\bibitem[Bonaventura and Ringler(2005)]{bona05Ay}
L.~Bonaventura and T.~D. Ringler.
\newblock Analysis of discrete shallow-water models on geodesic {Delaunay}
  grids with {C-Type} staggering.
\newblock \emph{Monthly Weather Review}, 133\penalty0 (8):\penalty0 2351--2373,
  2005.

\bibitem[Brecht et~al.(2019)Brecht, Bauer, Bihlo, Gay-Balmaz, and
  MacLachlan]{brecht2019variational}
R.~Brecht, W.~Bauer, A.~Bihlo, F.~Gay-Balmaz, and S.~MacLachlan.
\newblock Variational integrator for the rotating shallow-water equations on
  the sphere.
\newblock \emph{Quarterly Journal of the Royal Meteorological Society},
  145\penalty0 (720):\penalty0 1070--1088, 2019.

\bibitem[Brecht et~al.(2021)Brecht, Li, Bauer, and
  M{\'e}min]{brecht2021rotating}
R.~Brecht, L.~Li, W.~Bauer, and E.~M{\'e}min.
\newblock Rotating shallow water flow under location uncertainty with a
  structure-preserving discretization.
\newblock \emph{accepted manuscript, Journal of Advances in Modeling Earth
  Systems}, 2021.

\bibitem[Chen et~al.(2011)Chen, Gunzburger, and Ringler]{chen11Ay}
Q.~Chen, M.~Gunzburger, and T.~Ringler.
\newblock A scale-invariant formulation of the anticipated potential vorticity
  method.
\newblock \emph{Monthly Weather Review}, 139\penalty0 (8):\penalty0 2614--2629,
  2011.

\bibitem[Desbrun et~al.(2014)Desbrun, Gawlik, Gay-Balmaz, and
  Zeitlin]{desbrun2014variational}
M.~Desbrun, E.S. Gawlik, F.~Gay-Balmaz, and V.~Zeitlin.
\newblock Variational discretization for rotating stratified fluids.
\newblock \emph{Discrete \& Continuous Dynamical Systems-A}, 34\penalty0
  (2):\penalty0 477, 2014.

\bibitem[Flyer et~al.(2012)Flyer, Lehto, Blaise, Wright, and St-Cyr]{flye12a}
N.~Flyer, E.~Lehto, S.~Blaise, G.~B. Wright, and A.~St-Cyr.
\newblock A guide to {RBF}-generated finite differences for nonlinear
  transport: {S}hallow water simulations on a sphere.
\newblock \emph{Journal of Computational Physics}, 231\penalty0 (11):\penalty0
  4078--4095, 2012.

\bibitem[Gawlik and Gay-Balmaz(2020)]{GaGB2019}
E.S. Gawlik and F.~Gay-Balmaz.
\newblock A variational finite element discretization of compressible flow.
\newblock \emph{Foundations of Computational Mathematics}, 2020.

\bibitem[Gawlik et~al.(2011)Gawlik, Mullen, Pavlov, Marsden, and
  Desbrun]{gawlik2011geometric}
E.S. Gawlik, P.~Mullen, D.~Pavlov, J.E. Marsden, and M.~Desbrun.
\newblock Geometric, variational discretization of continuum theories.
\newblock \emph{Physica D: Nonlinear Phenomena}, 240\penalty0 (21):\penalty0
  1724--1760, 2011.

\bibitem[Gay-Balmaz and Holm(2014)]{GBHo2014}
F.~Gay-Balmaz and Darryl~D. Holm.
\newblock A geometric theory of selective decay with applications in {MHD}.
\newblock \emph{Nonlinearity}, 27:\penalty0 1747--1777, 2014.

\bibitem[Gay-Balmaz and Holm(2013)]{gay2013selective}
F.~Gay-Balmaz and D.D. Holm.
\newblock Selective decay by {C}asimir dissipation in inviscid fluids.
\newblock \emph{Nonlinearity}, 26\penalty0 (2):\penalty0 495, 2013.

\bibitem[Hairer et~al.(2006)Hairer, Lubich, and Wanner]{hair06Ay}
E.~Hairer, C.~Lubich, and G.~Wanner.
\newblock \emph{{Geometric numerical integration: structure-preserving
  algorithms for ordinary differential equations}}.
\newblock Springer, Berlin, 2006.

\bibitem[Heikes and Randall(1995{\natexlab{a}})]{heikes1995numericalP1}
R.~Heikes and D.~A Randall.
\newblock Numerical integration of the shallow-water equations on a twisted
  icosahedral grid. part i: Basic design and results of tests.
\newblock \emph{Monthly Weather Review}, 123\penalty0 (6):\penalty0 1862--1880,
  1995{\natexlab{a}}.

\bibitem[Heikes and Randall(1995{\natexlab{b}})]{heikes1995numerical}
R.~Heikes and D.A. Randall.
\newblock Numerical integration of the shallow-water equations on a twisted
  icosahedral grid. {Part II. A} detailed description of the grid and an
  analysis of numerical accuracy.
\newblock \emph{Monthly Weather Review}, 123\penalty0 (6):\penalty0 1881--1887,
  1995{\natexlab{b}}.

\bibitem[Holm et~al.(1998)Holm, Marsden, and Ratiu]{holm98Ay}
D.~D. Holm, J.~E. Marsden, and T.~S. Ratiu.
\newblock The {E}uler--{P}oincar\'{e} equations and semidirect products with
  applications to continuum theories.
\newblock \emph{Advances in Mathematics}, 137:\penalty0 1--81, 1998.

\bibitem[Kraichnan(1967)]{krai67Ay}
R.~H. Kraichnan.
\newblock Inertial ranges in two-dimensional turbulence.
\newblock \emph{The Physics of Fluids}, 10\penalty0 (7):\penalty0 1417--1423,
  1967.

\bibitem[Leimkuhler and Reich(2004)]{leim04Ay}
B.~Leimkuhler and S.~Reich.
\newblock \emph{{Simulating Hamiltonian dynamics}}.
\newblock Cambridge University Press, Cambridge, 2004.

\bibitem[Li et~al.(2020)Li, Wang, and Dong]{li2020analysis}
J.~Li, B.~Wang, and L.~Dong.
\newblock Analysis of and solution to the polar numerical noise within the
  shallow-water model on the latitude-longitude grid.
\newblock \emph{Journal of Advances in Modeling Earth Systems}, 12\penalty0
  (8):\penalty0 e2020MS002047, 2020.

\bibitem[Lilly(1971)]{lill71Ay}
D.~K. Lilly.
\newblock Numerical simulation of developing and decaying two-dimensional
  turbulence.
\newblock \emph{Journal of Fluid Mechanics}, 45\penalty0 (2):\penalty0
  395--415, 1971.

\bibitem[Marsden and West(2001)]{mars01a}
J.~E. Marsden and M.~West.
\newblock Discrete mechanics and variational integrators.
\newblock \emph{Acta Numerica 2001}, 10:\penalty0 357--514, 2001.

\bibitem[Matthaeus and Montgomery(1980)]{matthaeus1980selective}
W.H. Matthaeus and D.~Montgomery.
\newblock Selective decay hypothesis at high mechanical and magnetic reynolds
  numbers.
\newblock \emph{New York Academy of Sciences, Annals}, 357:\penalty0 203--222,
  1980.

\bibitem[McRae and Cotter(2014)]{mcra13Ay}
A.~T.~T. McRae and C.~J. Cotter.
\newblock Energy-and enstrophy-conserving schemes for the shallow-water
  equations, based on mimetic finite elements.
\newblock \emph{Quarterly Journal of the Royal Meteorological Society},
  140:\penalty0 2223--2234, 2014.

\bibitem[Nair et~al.(2021)Nair, Hanna, and Aureli]{nair2021selective}
A.G. Nair, J.~Hanna, and M.~Aureli.
\newblock Selective energy and enstrophy modification of two-dimensional
  decaying turbulence.
\newblock \emph{arXiv preprint arXiv:2108.01137}, 2021.

\bibitem[Natale and Cotter(2017)]{natale2017scale}
A.~Natale and C.~J. Cotter.
\newblock Scale-selective dissipation in energy-conserving finite-element
  schemes for two-dimensional turbulence.
\newblock \emph{Quarterly Journal of the Royal Meteorological Society},
  143\penalty0 (705):\penalty0 1734--1745, 2017.

\bibitem[Pavlov et~al.(2011)Pavlov, Mullen, Tong, Kanso, Marsden, and
  Desbrun]{pavlov2011structure}
D.~Pavlov, P.~Mullen, Y.~Tong, E.~Kanso, J.E. Marsden, and M.~Desbrun.
\newblock Structure-preserving discretization of incompressible fluids.
\newblock \emph{Physica D: Nonlinear Phenomena}, 240\penalty0 (6):\penalty0
  443--458, 2011.

\bibitem[Perot et~al.(2006)Perot, Vidovic, and Wesseling]{perot2006mimetic}
J.B. Perot, D.~Vidovic, and P.~Wesseling.
\newblock Mimetic reconstruction of vectors.
\newblock In \emph{Compatible Spatial Discretizations}, pages 173--188.
  Springer, 2006.

\bibitem[Ringler et~al.(2008)Ringler, Ju, and Gunzburger]{ring08Ay}
T.~Ringler, L.~Ju, and M.~Gunzburger.
\newblock {A multiresolution method for climate system modeling: application of
  spherical centroidal Voronoi tessellations}.
\newblock \emph{Ocean Dynamics}, 58:\penalty0 475--498, 2008.

\bibitem[Ringler and Randall(2002)]{ring02a}
T.~D. Ringler and D.~A. Randall.
\newblock A potential enstrophy and energy conserving numerical scheme for
  solution of the shallow-water equations on a geodesic grid.
\newblock \emph{Monthly Weather Review}, 130\penalty0 (5):\penalty0 1397--1410,
  2002.

\bibitem[Ringler et~al.(2010)Ringler, Thuburn, Klemp, and
  Skamarock]{ringler2010unified}
T.D. Ringler, J.~Thuburn, J.B. Klemp, and W.C. Skamarock.
\newblock A unified approach to energy conservation and potential vorticity
  dynamics for arbitrarily-structured c-grids.
\newblock \emph{Journal of Computational Physics}, 229\penalty0 (9):\penalty0
  3065--3090, 2010.

\bibitem[R{\'i}podas et~al.(2009)R{\'i}podas, Gassmann, F{\"o}rstner, Majewski,
  Giorgetta, Korn, Kornblueh, Wan, Z{\"a}ngl, Bonaventura, and
  Heinze]{Ripodas2009}
P.~R{\'i}podas, A.~Gassmann, J.~F{\"o}rstner, D.~Majewski, M.~Giorgetta,
  P.~Korn, L.~Kornblueh, H.~Wan, G.~Z{\"a}ngl, L.~Bonaventura, and T.~Heinze.
\newblock Icosahedral shallow water model ({ICOSWM}): results of shallow water
  test cases and sensitivity to model parameters.
\newblock \emph{Geoscientific Model Development}, 2\penalty0 (2):\penalty0
  231--251, 2009.

\bibitem[Sadourny and Basdevant(1985)]{sado85Ay}
R.~Sadourny and C.~Basdevant.
\newblock {Parameterization of subgrid scale barotropic and baroclinic eddies
  in quasi-geostrophic models: Anticipated potential vorticity method}.
\newblock \emph{Journal of Atmospheric Sciences}, 42\penalty0 (13):\penalty0
  1353--1363, 1985.

\bibitem[Shipton et~al.(2018)Shipton, Gibson, and Cotter]{shipton2018higher}
J.~Shipton, T.H. Gibson, and C.J. Cotter.
\newblock Higher-order compatible finite element schemes for the nonlinear
  rotating shallow water equations on the sphere.
\newblock \emph{Journal of Computational Physics}, 375:\penalty0 1121--1137,
  2018.

\bibitem[Shutts(2005)]{shutts2005kinetic}
G.~Shutts.
\newblock A kinetic energy backscatter algorithm for use in ensemble prediction
  systems.
\newblock \emph{Quarterly Journal of the Royal Meteorological Society: A
  journal of the atmospheric sciences, applied meteorology and physical
  oceanography}, 131\penalty0 (612):\penalty0 3079--3102, 2005.

\bibitem[Staniforth and Thuburn(2012)]{Staniforth12}
A.~Staniforth and J.~Thuburn.
\newblock Horizontal grids for global weather and climate prediction models: a
  review.
\newblock \emph{Quarterly Journal of the Royal Meteorological Society},
  138\penalty0 (662):\penalty0 1--26, 2012.

\bibitem[Thuburn et~al.(2014)Thuburn, Kent, and Wood]{thuburn2014cascades}
J.~Thuburn, J.~Kent, and N.~Wood.
\newblock Cascades, backscatter and conservation in numerical models of
  two-dimensional turbulence.
\newblock \emph{Quarterly Journal of the Royal Meteorological Society},
  140\penalty0 (679):\penalty0 626--638, 2014.

\bibitem[Wan and Nave(2016)]{wan16b}
A.~T.~S. Wan and J.-C. Nave.
\newblock On the arbitrarily long-term stability of conservative methods.
\newblock arXiv:1607.06160, 2016.

\bibitem[Warneford and Dellar(2014)]{warneford2014thermal}
E.S. Warneford and P.J. Dellar.
\newblock Thermal shallow water models of geostrophic turbulence in jovian
  atmospheres.
\newblock \emph{Physics of Fluids}, 26\penalty0 (1):\penalty0 016603, 2014.

\bibitem[Williamson et~al.(1992)Williamson, Drake, Hack, Jakob, and
  Swarztrauber]{wil92}
David~L. Williamson, John~B. Drake, James~J. Hack, R\"{u}diger Jakob, and
  Paul~N. Swarztrauber.
\newblock A standard test set for numerical approximations to the shallow water
  equations in spherical geometry.
\newblock \emph{Journal of Computational Physics}, 102\penalty0 (1):\penalty0
  211--224, 1992.

\bibitem[Wimmer et~al.(2020)Wimmer, Cotter, and Bauer]{wimmer2020energy}
G.A. Wimmer, C.J. Cotter, and W.~Bauer.
\newblock Energy conserving upwinded compatible finite element schemes for the
  rotating shallow water equations.
\newblock \emph{Journal of Computational Physics}, 401:\penalty0 109016, 2020.

\end{thebibliography}

\end{document}